\documentclass{amsart}

\usepackage{amssymb}
\usepackage{amstext}
\usepackage{amsmath}
\usepackage{amscd}
\usepackage{amsthm}
\usepackage{amsfonts}
\usepackage{enumerate}
\usepackage{graphicx}
\usepackage{latexsym}
\usepackage{mathrsfs}
\usepackage[all,color]{xy}
\usepackage{mathtools}
\xyoption{all}
\usepackage[usenames,dvipsnames]{xcolor}

\usepackage{pstricks}
\usepackage{lscape}
\usepackage{comment}

\newtheorem{theorem}{Theorem}[section]

\newtheorem{corollary}[theorem]{Corollary}

\newtheorem{proposition}[theorem]{Proposition}
\newtheorem{definition-proposition}[theorem]{Definition-Proposition}
\newtheorem{problem}[theorem]{Problem}

\theoremstyle{definition}
\newtheorem{definition}[theorem]{Definition}
\newtheorem{assumption}[theorem]{Assumption}

\newtheorem{example}[theorem]{Example}

\newcommand{\mm}{{\mathfrak{m}}}

\newcommand{\pp}{{\mathfrak{p}}}

\renewcommand{\AA}{\mathcal{A}}

\newcommand{\CC}{\mathcal{C}}
\newcommand{\CCC}{{\rm C}}
\newcommand{\DDD}{{\rm D}}
\newcommand{\KKK}{{\rm K}}
\newcommand{\OO}{{\mathcal O}}

\newcommand{\TT}{\mathcal{T}}

\newcommand{\UUU}{{\rm U}}

\newcommand{\FF}{\mathcal{F}}

\renewcommand{\c}{\vec{c}}

\newcommand{\x}{\vec{x}}
\newcommand{\y}{\vec{y}}
\newcommand{\z}{\vec{z}}
\newcommand{\w}{\vec{\omega}}

\newcommand{\A}{\mathbb{A}}
\newcommand{\D}{\mathbb{D}}
\newcommand{\E}{\mathbb{E}}
\renewcommand{\L}{\mathbb{L}}
\newcommand{\Z}{\mathbb{Z}}
\newcommand{\Q}{\mathbb{Q}}
\renewcommand{\P}{\mathbb{P}}

\newcommand{\SSS}{\mathbb{S}}

\newcommand{\X}{\mathbb{X}}

\newcommand{\bo}{\operatorname{b}\nolimits}
\newcommand{\sg}{\operatorname{sg}\nolimits}

\newcommand{\depth}{\operatorname{depth}\nolimits}

\newcommand{\soc}{\operatorname{soc}\nolimits}
\newcommand{\Ext}{\operatorname{Ext}\nolimits}

\newcommand{\Hom}{\operatorname{Hom}\nolimits}

\newcommand{\End}{\operatorname{End}\nolimits}

\newcommand{\gl}{\mathop{{\rm gl.dim}}\nolimits}

\newcommand{\op}{\operatorname{op}\nolimits}
\newcommand{\RHom}{\mathbf{R}\strut\kern-.2em\operatorname{Hom}\nolimits}
\newcommand{\RshHom}{\mathbf{R}\strut\kern-.2em\mathscr{H}\strut\kern-.3em\operatorname{om}\nolimits}
\newcommand{\shHom}{\mathscr{H}\strut\kern-.3em\operatorname{om}\nolimits}
\newcommand{\Lotimes}{\mathop{\stackrel{\mathbf{L}}{\otimes}}\nolimits}

\newcommand{\Spec}{\operatorname{Spec}\nolimits}

\newcommand{\ann}{\operatorname{ann}\nolimits}
\newcommand{\Proj}{\operatorname{Proj}\nolimits}

\newcommand{\SL}{\operatorname{SL}\nolimits}

\DeclareMathOperator{\moduleCategory}{{\rm mod}} \renewcommand{\mod}{\moduleCategory}
\DeclareMathOperator{\Mod}{{\rm Mod}}
\DeclareMathOperator{\proj}{{\rm proj}}

\DeclareMathOperator{\ind}{{\rm ind}}

\DeclareMathOperator{\coh}{{\rm coh}}
\DeclareMathOperator{\vect}{{\rm vect}}

\DeclareMathOperator{\qgr}{{\rm qgr}}
\DeclareMathOperator{\CM}{{\rm CM}}

\DeclareMathOperator{\add}{{\rm add}}

\newcommand{\cut}{\ar@{.}}
\newcommand{\rel}{\ar@{.}}

\numberwithin{equation}{section}

\SelectTips{cm}{}

\begin{document}
\title{Tilting Cohen-Macaulay representations}

\dedicatory{Dedicated to the memory of Ragnar-Olaf Buchweitz}

\author{Osamu Iyama}
\address{O. Iyama: Graduate School of Mathematics, Nagoya University, Chikusa-ku, Nagoya, 464-8602, Japan}
\email{iyama@math.nagoya-u.ac.jp}
\urladdr{http://www.math.nagoya-u.ac.jp/~iyama/}
\thanks{The author was partially supported by JSPS Grant-in-Aid for Scientific Research (B) 16H03923 and (S) 15H05738.}

\begin{abstract}
This is a survey on recent developments in Cohen-Macaulay representations via tilting and cluster tilting theory. 
We explain triangle equivalences between the singularity categories of Gorenstein rings and the derived
(or cluster) categories of finite dimensional algebras.
\end{abstract}
\maketitle

\section{Introduction}
The study of Cohen-Macaulay (CM) representations \cite{CR,Yo,Si,LW} is one of the active subjects in
representation theory and commutative algebra. It has fruitful connections to singularity theory,
algebraic geometry and physics.
This article is a survey on recent developments in this subject.

The first half of this article is spent for background materials, which were never written in one place.
In Section 2, we recall the notion of CM modules over Gorenstein rings,
and put them into the standard framework of triangulated categories.
This gives us powerful tools including Buchweitz's equivalence
between the stable category $\underline{\CM}R$ and the singularity category,
and Orlov's realization of the graded singularity category in the derived category, giving
a surprising connection between CM modules and algebraic geometry.
We also explain basic results including Auslander-Reiten duality stating that
$\underline{\CM}R$ is a Calabi-Yau triangulated category for a Gorenstein isolated singularity $R$,
and Gabriel's Theorem on quiver representations and its commutative counterpart
due to Buchweitz-Greuel-Schreyer.

In Section 3, we give a brief introduction to tilting and cluster tilting.
Tilting theory controls equivalences of derived categories, and played a central role in
Cohen-Macaulay approximation theory around 1990 \cite{ABu,AR2}.
The first main problem of this article is to find a tilting object in the stable category $\underline{\CM}^GR$ of a $G$-graded Gorenstein ring $R$.
This is equivalent to find a triangle equivalence
\begin{equation}\label{triangle equivalence1}
\underline{\CM}^GR\simeq\KKK^{\bo}(\proj\Lambda)
\end{equation}
with some ring $\Lambda$. It reveals a deep connection between rings $R$ and $\Lambda$.

The notion of $d$-cluster tilting was introduced in higher Auslander-Reiten theory. 
A Gorenstein ring $R$ is called $d$-CM-finite if there exists a $d$-cluster tilting object in $\CM R$.
This property is a natural generalization of CM-finiteness, and closely related to the existence of non-commutative crepant resolutions of Van den Bergh.
On the other hand, the $d$-cluster category $\CCC_d(\Lambda)$ of a finite dimensional algebra $\Lambda$ is a $d$-Calabi-Yau triangulated category containing a $d$-cluster tilting object, introduced in categorification of Fomin-Zelevinsky cluster algebras.
The second main problem of this article is to find a triangle equivalence
\begin{equation}\label{triangle equivalence2}
\underline{\CM} R\simeq\CCC_d(\Lambda)
\end{equation}
with some finite dimensional algebra $\Lambda$.
This implies that $R$ is $d$-CM-finite.

In the latter half of this article, we construct various triangle equivalences of the form
\eqref{triangle equivalence1} or \eqref{triangle equivalence2}.
In Section 4, we explain results in \cite{Ya} and \cite{BIY}. They assert that,
for a large class of $\Z$-graded Gorenstein rings $R$ in dimension 0 or 1,
there exist triangle equivalences \eqref{triangle equivalence1} for some algebras $\Lambda$.

There are no such general results in dimension greater than 1.
Therefore in the main Sections 5 and 6 of this article, we concentrate on special classes of
Gorenstein rings.
In Section 5, we explain results on Gorenstein rings obtained from classical and higher
preprojective algebras \cite{AIR,IO,Ki1,Ki2}.
A crucial observation is that certain Calabi-Yau algebras are higher preprojective algebras and higher Auslander algebras at the same time.
In Section 6, we explain results on CM modules on Geigle-Lenzing complete intersections and the derived categories of coherent sheaves on the associated stacks \cite{HIMO}.
They are higher dimensional generalizations of weighted projective lines of Geigle-Lenzing.

\section{Preliminaries}

\subsection{Notations}
We fix some conventions in this paper.
All modules are right modules. The composition of $f\colon X\to Y$ and $g\colon Y\to Z$ is denoted by $gf$.
For a ring $\Lambda$, we denote by $\mod\Lambda$ the category of finitely generated $\Lambda$-modules,
by $\proj\Lambda$ the category of finitely generated projective $\Lambda$-modules,
and by $\gl\Lambda$ the global dimension of $\Lambda$.
When $\Lambda$ is $G$-graded, we denote by $\mod^G\Lambda$ and $\proj^G\Lambda$
the $G$-graded version, whose morphisms are degree preserving.
We denote by $k$ an arbitrary field unless otherwise specified, and by $D$ the $k$-dual
or Matlis dual over a base commutative ring.

\subsection{Cohen-Macaulay modules}
We start with the classical notion of Cohen-Macaulay modules over commutative rings \cite{BrH,Ma}.

Let $R$ be a commutative noetherian ring.
The \emph{dimension} $\dim R$ of $R$ is the supremum of integers $n\ge0$ such that 
there exists a chain $\pp_0\subsetneq\pp_1\subsetneq\cdots\subsetneq\pp_n$ of prime ideals of $R$.
The \emph{dimension} $\dim M$ of $M\in\mod R$ is the dimension $\dim(R/\ann M)$ of
the factor ring $R/\ann M$, where $\ann M$ is the annihilator of $M$.

The notion of depth is defined locally.
Assume that $R$ is a local ring with maximal ideal $\mm$ and $M\in\mod R$ is non-zero.
An element $r\in\mm$ is called \emph{$M$-regular} if the multiplication map
$r\colon M\to M$ is injective. A sequence $r_1,\ldots,r_n$ of elements in $\mm$ is called
an \emph{$M$-regular sequence} of length $n$ if $r_i$ is $(M/(r_1,\ldots,r_{i-1})M)$-regular for all $1\le i\le n$.
The \emph{depth} $\depth M$ of $M$ is the supremum of the length of $M$-regular sequences.
This is  given by the simple formula
\[\depth M=\inf\{i\ge0\mid \Ext^i_R(R/\mm,M)\neq0\}.\]
The inequalities $\depth M\le\dim M\le\dim R$ hold. We call $M$ (maximal) \emph{Cohen-Macaulay}
(or \emph{CM}) if the equality $\depth M=\dim R$ holds or $M=0$.

When $R$ is not necessarily local, $M\in\mod R$ is called \emph{CM} if $M_{\mm}$ is a
CM $R_{\mm}$-module for all maximal ideals $\mm$ of $R$.
The ring $R$ is called \emph{CM} if it is CM as an $R$-module.
The ring $R$ is called \emph{Gorenstein} (resp.\ \emph{regular}) if $R_{\mm}$ has
finite injective dimension as an $R_{\mm}$-module (resp.\ $\gl R_{\mm}<\infty$)
for all maximal ideals $\mm$ of $R$.
In this case, the injective (resp.\ global) dimension coincides with $\dim R_{\mm}$,
but this is not true in the more general setting below. The following hierarchy is basic.
\[\xymatrix{\mbox{Regular rings}\ar@{=>}[r]&\mbox{Gorenstein rings}\ar@{=>}[r]&\mbox{Cohen-Macaulay rings}}\]
We will study CM modules over Gorenstein rings. Since we apply methods in representation theory,
it is more reasonable to work in the following wider framework.

\begin{definition}\cite{Iw,EJ}\label{define IG}
Let $\Lambda$ be a (not necessarily commutative) noetherian ring, and $d\ge0$ an integer.
We call $\Lambda$ ($d$-)\emph{Iwanaga-Gorenstein} (or \emph{Gorenstein}) if $\Lambda$ has injective 
dimension at most $d$ as a $\Lambda$-module, and also as a $\Lambda^{\op}$-module.
\end{definition}

Clearly, a commutative noetherian ring $R$ is Iwanaga-Gorenstein if and only if it is Gorenstein and $\dim R<\infty$. Note that there are various definitions of non-commutative
Gorenstein rings, e.g.\ \cite{AS,CR,FGR,GN,IW}.
Although Definition \ref{define IG} is much weaker than them, it is sufficient for the aim of this paper.

Noetherian rings with finite global dimension are analogues of regular rings, and form special classes of Iwanaga-Gorenstein rings.
The first class consists of \emph{semisimple rings}
(i.e.\ rings $\Lambda$ with $\gl\Lambda=0$), which are products of matrix rings over division rings
by the Artin-Wedderburn Theorem.
The next class consists of \emph{hereditary rings} (i.e.\ rings $\Lambda$ with $\gl\Lambda\le1$),
which are obtained from quivers.

\begin{definition}\cite{ASS}
A \emph{quiver} is a quadruple $Q=(Q_0,Q_1,s,t)$ consisting of sets $Q_0$, $Q_1$ and
maps $s,t\colon Q_1\to Q_0$. We regard each element in $Q_0$ as a vertex, and $a\in Q_1$
as an arrow with source $s(a)$ and target $t(a)$.
A \emph{path} of length 0 is a vertex, and a \emph{path} of length $\ell(\ge1)$ is a sequence
$a_1a_2\cdots a_\ell$ of arrows satisfying $t(a_i)=s(a_{i+1})$ for each $1\le i<\ell$.

For a field $k$, the \emph{path algebra} $kQ$ is defined as follows:
It is a $k$-vector space with basis consisting of all paths on $Q$. For paths $p=a_1\cdots a_\ell$
and $q=b_1\cdots b_m$, we define $pq=a_1\cdots a_\ell b_1\cdots b_m$ if $t(a_\ell)=s(b_1)$,
and $pq=0$ otherwise.
\end{definition}

Clearly $\dim_k(kQ)$ is finite if and only if $Q$ is \emph{acyclic} (that is, there are no paths $p$
of positive length satisfying $s(p)=t(p)$).

\begin{example}\label{example of IG}
\begin{enumerate}[\rm(a)]
\item \cite{ASS} The path algebra $kQ$ of a finite quiver $Q$ is hereditary.
Conversely, any finite dimensional hereditary algebra over an algebraically closed field $k$
is Morita equivalent to $kQ$ for some acyclic quiver $Q$.
\item A finite dimensional $k$-algebra $\Lambda$ is $0$-Iwanaga-Gorenstein if and only if
$\Lambda$ is \emph{self-injective}, that is, $D\Lambda$ is projective as a $\Lambda$-module,
or equivalently, as a $\Lambda^{\op}$-module.
For example, the group ring $kG$ of a finite group $G$ is self-injective.
\item \cite{IW,CR} Let $R$ be a CM local ring with canonical module $\omega$ and dimension $d$.
An $R$-algebra $\Lambda$ is called an \emph{$R$-order} if it is CM as an $R$-module.
Then an $R$-order $\Lambda$ is $d$-Iwanaga-Gorenstein if and only if $\Lambda$ is
a \emph{Gorenstein order}, i.e.\ $\Hom_R(\Lambda,\omega)$ is
projective as a $\Lambda$-module, or equivalently, as a $\Lambda^{\op}$-module.

An $R$-order $\Lambda$ is called \emph{non-singular} if $\gl\Lambda=d$.
They are classical objects for the case $d=0,1$ \cite{CR}, and studied for $d=2$  \cite{RV1}.
Non-singular orders are closely related to cluster tilting explained in Section \ref{sec: Cluster tilting theory}.
\end{enumerate}
\end{example}

\subsection{The triangulated category of Cohen-Macaulay modules}

CM modules can be defined naturally also for Iwanaga-Gorenstein rings.

\begin{definition}\label{define CM}
Let $\Lambda$ be an Iwanaga-Gorenstein ring.
We call $M\in\mod\Lambda$ (maximal) \emph{Cohen-Macaulay} (or \emph{CM})
if $\Ext^i_\Lambda(M,\Lambda)=0$ holds for all $i>0$.
We denote by $\CM\Lambda$ the category of CM $\Lambda$-modules.
\end{definition}

We also deal with graded rings and modules.
For an abelian group $G$ and a $G$-graded Iwanaga-Gorenstein ring $\Lambda$,
we denote by $\CM^G\Lambda$ the full subcategory of $\mod^G\Lambda$ consisting of all $X$
which belong to $\CM\Lambda$ as ungraded $\Lambda$-modules.

When $\Lambda$ is commutative Gorenstein, Definition \ref{define CM} is one of the well-known equivalent
conditions of CM modules.
Note that, in a context of Gorenstein homological algebra \cite{ABr,EJ}, CM modules are also called
Gorenstein projective, Gorenstein dimension zero, or totally reflexive.

\begin{example}
\begin{enumerate}[\rm(a)]
\item Let $\Lambda$ be a noetherian ring with $\gl\Lambda<\infty$.
Then $\CM\Lambda=\proj\Lambda$.
\item Let $\Lambda$ be a finite dimensional self-injective $k$-algebra. Then $\CM\Lambda=\mod\Lambda$.
\item Let $\Lambda$ be a Gorenstein $R$-order in Example \ref{example of IG}(c).
Then CM $\Lambda$-modules are precisely $\Lambda$-modules that are CM as $R$-modules.
\end{enumerate}
\end{example}

We study the category $\CM^G\Lambda$ from the point of view of triangulated categories.
We start with Quillen's exact categories (see \cite{Bue} for a more axiomatic definition).

\begin{definition}\cite{H}
\begin{enumerate}[\rm(a)]
\item An \emph{exact category} is a full subcategory $\FF$ of an abelian category $\AA$
such that, for each exact sequence $0\to X\to Y\to Z\to0$ in $\AA$ with $X,Z\in\FF$, we have $Y\in\FF$.
In this case, we say that $X\in \FF$ is \emph{projective}
if $\Ext^1_{\AA}(X,\FF)=0$ holds. Similarly we define \emph{injective} objects in $\FF$.
\item An exact category $\FF$ in $\AA$ is called \emph{Frobenius} if:
\begin{enumerate}[$\bullet$]
\item an object in $\FF$ is projective if and only if it is injective,
\item any $X\in\FF$ admits exact sequences $0\to Y\to P\to X\to 0$ and $0\to X\to I\to Z\to 0$
in $\AA$ such that $P$ and $I$ are projective in $\FF$ and $Y,Z\in\FF$.
\end{enumerate}
\item The \emph{stable category} $\underline{\FF}$ has the same objects as $\FF$, and the morphisms are 
given by $\underline{\Hom}_{\FF}(X,Y)=\Hom_{\FF}(X,Y)/P(X,Y)$,
where $P(X,Y)$ is the subgroup consisting of morphisms which factor through projective objects in $\FF$.
\end{enumerate}
\end{definition}

Frobenius categories are ubiquitous in algebra. Here we give two examples.

\begin{example}
\begin{enumerate}[\rm(a)]
\item For a $G$-graded Iwanaga-Gorenstein ring $\Lambda$, the category $\CM^G\Lambda$ of
$G$-graded Cohen-Macaulay $\Lambda$-modules is a Frobenius category.
\item For an additive category $\AA$, the category $\CCC(\AA)$ of chain complexes in $\AA$
is a Frobenius category, whose stable category is the homotopy category $\KKK(\AA)$.
\end{enumerate}
\end{example}

A \emph{triangulated category} is a triple of an additive category $\TT$,
an autoequivalence $[1]\colon \TT\to\TT$ (called \emph{suspension}) and a class of diagrams
$X\xrightarrow{f}Y\xrightarrow{g}Z\xrightarrow{h}X[1]$ (called \emph{triangles})
satisfying a certain set of axioms.
There are natural notions of functors and equivalences between triangulated categories, called
\emph{triangle functors} and \emph{triangle equivalences}. For details, see e.g.\ \cite{H,N}.
Typical examples of triangulated categories are given by the homotopy category $\KKK(\AA)$
of an additive category $\AA$ and the derived category $\DDD(\AA)$ of an abelian category $\AA$.

A standard construction of triangulated categories is given by the following.

\begin{theorem}\cite{H}\label{happel algebraic}
The stable category $\underline{\FF}$ of a Frobenius category $\FF$ has a canonical structure
of a triangulated category.
\end{theorem}

Such a triangulated category is called \emph{algebraic}.
Note that the suspension functor $[1]$ of $\underline{\FF}$ is given by the cosyzygy.
Thus the $i$-th suspension $[i]$ is the $i$-th cosyzygy for $i\ge0$, and the $(-i)$-th syzygy for $i<0$.
We omit other details.

As a summary, we obtain the following.

\begin{corollary}
Let $G$ be an abelian group and $\Lambda$ a $G$-graded Iwanaga-Gorenstein ring.
Then $\CM^G\Lambda$ is a Frobenius category, and therefore the stable category
$\underline{\CM}^G\Lambda$ has a canonical structure of a triangulated category.
\end{corollary}

We denote by $\DDD^{\bo}(\mod^G\Lambda)$ the bounded derived category of $\mod ^G\Lambda$,
and by $\KKK^{\bo}(\proj^G\Lambda)$ the bounded homotopy category of $\proj^G\Lambda$.
We regard $\KKK^{\bo}(\proj^G\Lambda)$ as a thick subcategory of $\DDD^{\bo}(\mod^G\Lambda)$.
The \emph{stable derived category} \cite{Bu} or the \emph{singularity category} \cite{O} is defined as
the Verdier quotient
\[\DDD_{\sg}^G(\Lambda)=\DDD^{\bo}(\mod^G\Lambda)/\KKK^{\bo}(\proj^G\Lambda).\]
This is enhanced by the Frobenius category $\CM^G\Lambda$ as the following result shows.

\begin{theorem} \cite{Bu,Ric2,KV2}\label{stable is singularity}
Let $G$ be an abelian group and $\Lambda$ a $G$-graded Iwanaga-Gorenstein ring.
Then there is a triangle equivalence $\DDD_{\sg}^G(\Lambda)\simeq\underline{\CM}^G\Lambda$.
\end{theorem}

Let us recall the following notion \cite{BrH}.

\begin{definition}\label{define a-invariant}
Let $G$ be an abelian group and $R$ a $G$-graded Gorenstein ring with $\dim R=d$
such that $R_0=k$ is a field and $\bigoplus_{i\neq0}R_i$ is an ideal of $R$. The \emph{$a$-invariant} $a\in G$
(or \emph{Gorenstein parameter} $-a\in G$) is an element satisfying
$\Ext^d_R(k,R(a))\simeq k$ in $\mod^{\Z}R$.
\end{definition}

For a $G$-graded noetherian ring $\Lambda$, let
\begin{equation}\label{define qgr}
\qgr\Lambda=\mod^{G}\Lambda/\mod^{G}_0\Lambda
\end{equation}
be the Serre quotient of $\mod^{G}\Lambda$ by the subcategory $\mod^{G}_0\Lambda$ of $G$-graded $\Lambda$-modules of finite length \cite{AZ}.
This is classical in projective geometry. In fact, for a $\Z$-graded commutative noetherian ring $R$ generated in degree 1, $\qgr R$ is the category $\coh X$ of coherent sheaves on the scheme $X=\Proj R$ \cite{Se}.

The following result realizes $\DDD_{\sg}^{\Z}(R)$ and
$\DDD^{\bo}(\qgr R)$ inside of $\DDD^{\bo}(\mod^{\Z}R)$,
where $\mod^{\ge n}R$ is the full subcategory
of $\mod^{\Z}R$ consisting of all $X$ satisfying $X=\bigoplus_{i\ge n}X_i$, and $(-)^*$ is the duality
$\RHom_R(-,R)\colon \DDD^{\bo}(\mod^{\Z}R)\to\DDD^{\bo}(\mod^{\Z}R)$. 

\begin{theorem}\cite{O,IYa}\label{orlov embedding}
Let $R=\bigoplus_{i\ge0}R_i$ be a $\Z$-graded Gorenstein ring such that $R_0$ is a field,
and $a$ the $a$-invariant of $R$.
\begin{enumerate}[\rm(a)]
\item There is a triangle equivalence $\DDD^{\bo}(\mod^{\ge0}R)\cap\DDD^{\bo}(\mod^{\ge1}R)^*\simeq\DDD_{\sg}^{\Z}(R)$.
\item There is a triangle equivalence $\DDD^{\bo}(\mod^{\ge0}R)\cap\DDD^{\bo}(\mod^{\ge a+1}R)^*\simeq\DDD^{\bo}(\qgr R)$.
\end{enumerate}
\end{theorem}

Therefore if $a=0$, then $\DDD_{\sg}^{\Z}(R)\simeq\DDD^{\bo}(\qgr R)$. If $a<0$ (resp. $a>0$),
then there is a fully faithful triangle functor $\DDD_{\sg}^{\Z}(R)\to\DDD^{\bo}(\qgr R)$
(resp. $\DDD^{\bo}(\qgr R)\to\DDD_{\sg}^{\Z}(R)$).
This gives a new connection between CM representations and algebraic geometry.

\subsection{Representation theory}\label{section: quiver representations}

We start with Auslander-Reiten theory.

Let $R$ be a commutative ring, and $D$ the Matlis duality.
A triangulated category $\TT$ is called \emph{$R$-linear} if
each morphism set $\Hom_{\TT}(X,Y)$ has an $R$-module structure and the composition
$\Hom_{\TT}(X,Y)\times\Hom_{\TT}(Y,Z)\to\Hom_{\TT}(X,Z)$ is $R$-bilinear. It is called
\emph{Hom-finite} if each morphism set has finite length as an $R$-module.

\begin{definition}\cite{RV2}
A \emph{Serre functor} is an $R$-linear autoequivalence $\SSS\colon \TT\to\TT$ such that there exists
a functorial isomorphism $\Hom_{\TT}(X,Y)\simeq D\Hom_{\TT}(Y,\SSS X)$ for any $X,Y\in\TT$
(called \emph{Auslander-Reiten duality} or \emph{Serre duality}).
The composition $\tau=\SSS\circ[-1]$ is called the \emph{AR translation}.

For $d\in\Z$, we say that $\TT$ is \emph{$d$-Calabi-Yau} if $[d]$ gives a Serre functor of $\TT$.
\end{definition}

A typical example of a Serre functor is given by a smooth projective variety $X$ over a field $k$. 
In this case, $\DDD^{\bo}(\coh X)$ has a Serre functor $-\otimes_X\omega[d]$, where
$\omega$ is the canonical bundle of $X$ and $d$ is the dimension of $X$ \cite{Hu}.

\begin{example}\cite{H,BIY}\label{serre functor of derived category}
Let $\Lambda$ be a finite dimensional $k$-algebra.
Then $\KKK^{\bo}(\proj\Lambda)$ has a Serre functor if and only if $\Lambda$ is Iwanaga-Gorenstein,
and $\DDD^{\bo}(\mod\Lambda)$ has a Serre functor if and only if $\gl\Lambda<\infty$.
In both cases, the Serre functor is given by $\nu=-\Lotimes_\Lambda(D\Lambda)$, and
the AR translation is given by $\tau=\nu\circ[-1]$.
\end{example}

For AR theory of CM modules, we need the following notion.

\begin{definition}\label{define CM_0}
Let $R$ be a Gorenstein ring with $\dim R=d$. We denote by $\CM_0R$ the full subcategory of $\CM R$
consisting of all $X$ such that $X_{\pp}\in\proj R_{\pp}$ holds for all $\pp\in\Spec R$ with $\dim R_{\pp}<d$.
When $R$ is local, such an $X$ is called \emph{locally free on the punctured spectrum} \cite{Yo}.
If $R$ is $G$-graded, we denote by $\CM_0^GR$ the full subcategory of $\CM^GR$ consisting of
all $X$ which belong to $\CM_0R$ as ungraded $R$-modules.
\end{definition}

As before, $\CM^G_0R$ is a Frobenius category, and $\underline{\CM}^G_0R$ is a triangulated category.
Note that $\CM_0R=\CM R$ holds if and only if $R$ satisfies Serre's (R$_{d-1}$) condition
(i.e.\ $R_{\pp}$ is regular for all $\pp\in\Spec R$ with $\dim R_{\pp}<d$).
This means that $R$ has at worst an isolated singularity if $R$ is local.

The following is a fundamental theorem of CM representations.

\begin{theorem}\cite{Au1,AR1}\label{AR duality for CM}
Let $R$ be a Gorenstein ring with $\dim R=d$. Then $\underline{\CM}_0R$ is
a $(d-1)$-Calabi-Yau triangulated category. If $R$ is $G$-graded and has
an $a$-invariant $a\in G$, then $\underline{\CM}_0^GR$ has a Serre functor $(a)[d-1]$.
\end{theorem}

Let us introduce a key notion in Auslander-Reiten theory.
We call an additive category $\CC$ \emph{Krull-Schmidt} if any object in $\CC$ is isomorphic to 
a finite direct sum of objects whose endomorphism rings are local.
We denote by $\ind\CC$ the set of isomorphism classes of indecomposable objects in $\CC$.

\begin{definition}\cite{H}
Let $\TT$ be a Krull-Schmidt triangulated category.
We call a triangle $X\xrightarrow{f} Y\xrightarrow{g} Z\xrightarrow{h} X[1]$ in $\TT$
an \emph{almost split triangle} if:
\begin{enumerate}[$\bullet$]
\item $X$ and $Z$ are indecomposable, and $h\neq0$ (i.e.\ the triangle does not split).
\item Any morphism $W\to Z$ which is not a split epimorphism factors through $g$.
\item Any morphism $X\to W$ which is not a split monomorphism factors through $f$.
\end{enumerate}
We say that \emph{$\TT$ has almost split triangles} if for any indecomposable object $X$
(resp.\ $Z$), there is an almost split triangle $X\to Y\to Z\to X[1]$.
\end{definition}

There is a close connection between almost split triangles and Serre functors.

\begin{theorem}\cite{RV2}
Let $\TT$ be an $R$-linear Hom-finite Krull-Schmidt triangulated category.
Then $\TT$ has a Serre functor if and only if $\TT$ has almost split triangles. In this case,
$X\simeq\tau Z$ holds in each almost split triangle $X\to Y\to Z\to X[1]$ in $\TT$.
\end{theorem}

When $\TT$ has almost split triangles, one can define the \emph{AR quiver} of $\TT$, which has $\ind\TT$ as the set of vertices. It describes the structure of $\TT$ (see \cite{H}).
Similarly, \emph{almost split sequences} and the AR quiver
are defined for exact categories \cite{ASS,LW}.

In the rest of this section, we discuss the following notion.

\begin{definition}
A finite dimensional $k$-algebra $\Lambda$ is called \emph{representation-finite} if
$\ind(\mod\Lambda)$ is a finite set. It is also said to be \emph{of finite representation type}.
\end{definition}

The classification of representation-finite algebras was one of the main subjects in the 1980s.
Here we recall only one theorem, and refer to \cite{GR} for further results.

A \emph{Dynkin quiver} (resp. \emph{extended Dynkin quiver}) is a quiver obtained by
orienting each edge of one of the following diagrams $A_n$, $D_n$ and $E_n$
(resp. $\widetilde{A}_n$, $\widetilde{D}_n$ and $\widetilde{E}_n$).
{\small
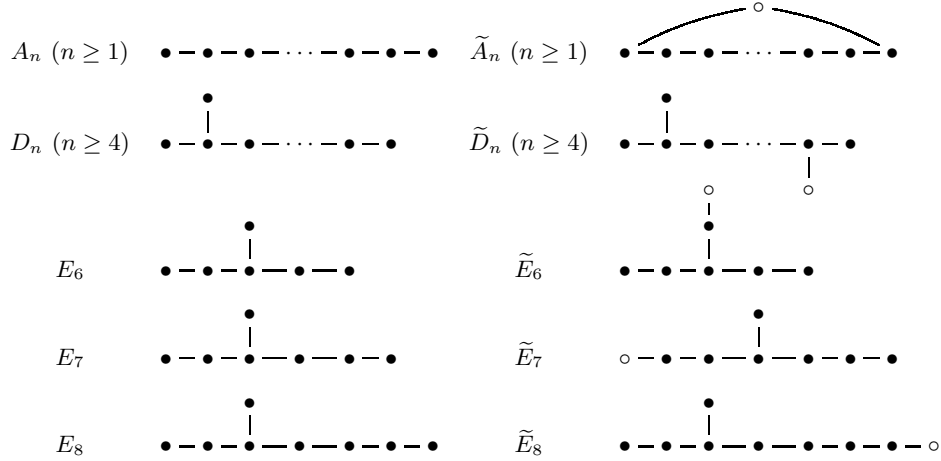
\begin{equation}\label{dynkin}
\xymatrix@C0.2cm@R0.15cm{
&&&&&&&&
&&&&\circ\ar@/^1mm/@{-}[drrr]\ar@/^-1mm/@{-}[dlll]&&&&\\
A_n\ (n\ge1)&\bullet\ar@{-}[r]&\bullet\ar@{-}[r]&\bullet\ar@{-}[r]&\cdots\ar@{-}[r]&\bullet\ar@{-}[r]&\bullet\ar@{-}[r]&\bullet&
\widetilde{A}_n\ (n\ge1)&\bullet\ar@{-}[r]&\bullet\ar@{-}[r]&\bullet\ar@{-}[r]&\cdots\ar@{-}[r]&\bullet\ar@{-}[r]&\bullet\ar@{-}[r]&\bullet\\
&&\bullet&&&&&&
&&\bullet\\
D_n\ (n\ge4)&\bullet\ar@{-}[r]&\bullet\ar@{-}[r]\ar@{-}[u]&\bullet\ar@{-}[r]&\cdots\ar@{-}[r]&\bullet\ar@{-}[r]&\bullet&&
\widetilde{D}_n\ (n\ge4)&\bullet\ar@{-}[r]&\bullet\ar@{-}[r]\ar@{-}[u]&\bullet\ar@{-}[r]&\cdots\ar@{-}[r]&\bullet\ar@{-}[r]\ar@{-}[d]&\bullet\\
&&&&&&&&
&&&\circ&&\circ\\
&&&\bullet&&&&&
&&&\bullet\ar@{-}[u]\\
E_6&\bullet\ar@{-}[r]&\bullet\ar@{-}[r]&\bullet\ar@{-}[r]\ar@{-}[u]&\bullet\ar@{-}[r]&\bullet&&&
\widetilde{E}_6&\bullet\ar@{-}[r]&\bullet\ar@{-}[r]&\bullet\ar@{-}[r]\ar@{-}[u]&\bullet\ar@{-}[r]&\bullet\\
&&&\bullet&&&&&
&&&&\bullet\\
E_7&\bullet\ar@{-}[r]&\bullet\ar@{-}[r]&\bullet\ar@{-}[r]\ar@{-}[u]&\bullet\ar@{-}[r]&\bullet\ar@{-}[r]&\bullet&&
\widetilde{E}_7&\circ\ar@{-}[r]&\bullet\ar@{-}[r]&\bullet\ar@{-}[r]&\bullet\ar@{-}[r]\ar@{-}[u]&\bullet\ar@{-}[r]&\bullet\ar@{-}[r]&\bullet\\
&&&\bullet&&&&&
&&&\bullet\\
E_8&\bullet\ar@{-}[r]&\bullet\ar@{-}[r]&\bullet\ar@{-}[r]\ar@{-}[u]&\bullet\ar@{-}[r]&\bullet\ar@{-}[r]&\bullet\ar@{-}[r]&\bullet&
\widetilde{E}_8&\bullet\ar@{-}[r]&\bullet\ar@{-}[r]&\bullet\ar@{-}[r]\ar@{-}[u]&\bullet\ar@{-}[r]&\bullet\ar@{-}[r]&\bullet\ar@{-}[r]&\bullet\ar@{-}[r]&\circ}
\end{equation}}
Now we are able to state Gabriel's Theorem below. 
For results in non-Dynkin case, we refer to Kac's theorem in \cite{GR}.

\begin{theorem}\cite{ASS}\label{gabriel}
Let $Q$ be a connected acyclic quiver and $k$ a field. Then $kQ$ is representation-finite
if and only if $Q$ is Dynkin. In this case, there is a bijection between $\ind(\mod kQ)$ and the set
$\Phi_+$ of positive roots in the root system of $Q$.
\end{theorem}

For a quiver $Q$, we define a new quiver $\Z Q$: The set of vertices is $\Z\times Q_0$.
The arrows are $(\ell,a)\colon (\ell,s(a))\to(\ell,t(a))$ and $(\ell,a^*)\colon (\ell,t(a))\to(\ell+1,s(a))$
for each $\ell\in\Z$ and $a\in Q_1$.
For example, if $Q=[1\xrightarrow{a} 2\xrightarrow{b} 3]$, then $\Z Q$ is as follows.
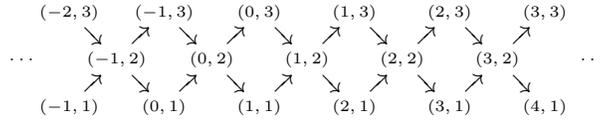
\[{\tiny\begin{xy}
0;<6pt,0pt>:<0pt,6pt>:: 
(-9,3) *+{\cdots},
(-6,6) *+{(-2,3)} ="-23",
(-6,0) *+{(-1,1)} ="-11",
(-3,3) *+{(-1,2)} ="-12",
(0,6) *+{(-1,3)} ="-13",
(0,0) *+{(0,1)} ="01",
(3,3) *+{(0,2)} ="02",
(6,6) *+{(0,3)} ="03",
(6,0) *+{(1,1)} ="11",
(9,3) *+{(1,2)} ="12",
(12,6) *+{(1,3)} ="13",
(12,0) *+{(2,1)} ="21",
(15,3) *+{(2,2)} ="22",
(18,6) *+{(2,3)} ="23",
(18,0) *+{(3,1)} ="31",
(21,3) *+{(3,2)} ="32",
(24,6) *+{(3,3)} ="33",
(24,0) *+{(4,1)} ="41",
(27,3) *+{\cdots},
\ar "-11";"-12"
\ar "-12";"-13"
\ar "01";"02"
\ar "02";"03"
\ar "11";"12"
\ar "12";"13"
\ar "21";"22"
\ar "22";"23"
\ar "31";"32"
\ar "32";"33"
\ar "-23";"-12"
\ar "-12";"01"
\ar "-13";"02"
\ar "02";"11"
\ar "03";"12"
\ar "12";"21"
\ar "13";"22"
\ar "22";"31"
\ar "23";"32"
\ar "32";"41"
\end{xy}}\]
If the underlying graph $\Delta$ of $Q$ is a tree, then $\Z Q$ depends only on $\Delta$.
Thus $\Z Q$ is written as $\Z\Delta$.

The AR quiver of $\DDD^{\bo}(\mod kQ)$ has a simple description \cite{H}.

\begin{proposition}\cite{H}\label{happel ZQ}
\begin{enumerate}[\rm(a)]
\item Let $\Lambda$ be a finite dimensional hereditary algebra. Then there is a bijection
$\ind(\mod\Lambda)\times\Z\to\ind\DDD^{\bo}(\mod\Lambda)$ given by $(X,i)\mapsto X[i]$.
\item For each Dynkin quiver $Q$, the AR quiver of $\DDD^{\bo}(\mod kQ)$ is $\Z Q^{\op}$. 
Moreover, the category $\DDD^{\bo}(\mod kQ)$ is presented by the quiver $\Z Q^{\op}$ with \emph{mesh relations}.
\end{enumerate}
\end{proposition}

Note that $\Z Q$ has an automorphism $\tau$ given by $\tau(\ell,i)=(\ell-1,i)$ for $(\ell,i)\in\Z\times Q_0$,
which corresponds to the AR translation.

Now we discuss CM-finiteness.
For an additive category $\CC$ and an object $M\in\CC$, we denote by $\add M$
the full subcategory of $\CC$ consisting of direct summands of finite direct sum of copies of $M$.
We call $M$ an \emph{additive generator} of $\CC$ if $\CC=\add M$.

\begin{definition}\label{define CM-finite}
An Iwanaga-Gorenstein ring $\Lambda$ is called \emph{CM-finite} if $\CM\Lambda$ has an additive generator $M$. In this case, we call $\End_\Lambda(M)$ the \emph{Auslander algebra}.
When $\CM\Lambda$ is Krull-Schmidt, $\Lambda$ is CM-finite if and only if $\ind(\CM\Lambda)$ is a finite set.
It is also said to be \emph{of finite CM type} or \emph{representation-finite}.
\end{definition}

Let us recall the classification of CM-finite Gorenstein rings given in the 1980s.
Let $k$ be an algebraically closed field of characteristic zero.
A hypersurface $R=k[[x,y,z_2\ldots,z_d]]/(f)$ is called a \emph{simple singularity} if 
\begin{equation}\label{simple equation}
f=\left\{\begin{array}{ll}
x^{n+1}+y^2+z_2^2+\cdots+z_d^2&A_n\\
x^{n-1}+xy^2+z_2^2+\cdots+z_d^2&D_n\\
x^4+y^3+z_2^2+\cdots+z_d^2&E_6\\
x^3y+y^3+z_2^2+\cdots+z_d^2&E_7\\
x^5+y^3+z_2^2+\cdots+z_d^2&E_8.
\end{array}\right.
\end{equation}
We are able to state the following result. We refer to \cite{LW} for results in positive characteristic.

\begin{theorem}\cite{BGS,Kn}\label{bgs}
Let $R$ be a complete local Gorenstein ring containing the residue field $k$,
which is an algebraically closed field of characteristic zero.
Then $R$ is CM-finite if and only if it is a simple singularity.
\end{theorem}

We will see that tilting theory explains why Dynkin quivers appear in both
Theorems \ref{gabriel} and \ref{bgs} (see Example \ref{simple curve singularity}
and Corollary \ref{simple surface singularity} below).

Now we describe the AR quivers of simple singularities.
Recall that each quiver $Q$ gives a new quiver $\Z Q$.
For an automorphism $\phi$ of $\Z Q$, an orbit quiver $\Z Q/\phi$ is naturally defined.
For example, $\Z Q/\tau$ is the \emph{double} $\overline{Q}$ of $Q$ obtained by
adding an inverse arrow $a^*\colon j\to i$ for each arrow $a\colon i\to j$. 

\begin{proposition}\cite{Yo,DW}\label{AR quiver of simple singularity}
Let $R$ be a simple singularity with $\dim R=d$. Then the AR quiver of $\underline{\CM} R$ is $\Z\Delta/\phi$,
where $\Delta$ and $\phi$ are given as follows.
\begin{enumerate}[\rm(a)]
\item If $d$ is even, then $\Delta$ is the Dynkin diagram of the same type as $R$, and $\phi=\tau$.
\item If $d$ is odd, then $\Delta$ and $\phi$ are given as follows.
{\small\[\begin{array}{|c||c|c|c|c|c|c|c|}
\hline
R&A_{2n-1}&A_{2n}&D_{2n}&D_{2n+1}&E_6&E_7&E_8\\ \hline\hline
\Delta&D_{n+1}&A_{2n}&D_{2n}&A_{4n-1}&E_6&E_7&E_8\\ \hline
\phi&\tau\iota&\tau^{1/2}&\tau^2&\tau\iota&\tau\iota&\tau^2&\tau^2\\ \hline
\end{array}\]}
Here $\iota$ is the involution of $\Z\Delta$ induced by the non-trivial involution of $\Delta$,
and $\tau^{1/2}$ is the automorphism of $\Z A_{2n}$ satisfying $(\tau^{1/2})^2=\tau$.
\end{enumerate}
\end{proposition}

In dimension 2, simple singularities (over a sufficiently large field) have
an alternative description as invariant subrings.
This enables us to draw the AR quiver of the category $\CM R$ systematically.
  
\begin{example}\label{quotient singularity2}\cite{Au2,LW}
Let $k[[u,v]]$ be a formal power series ring over a field $k$ and $G$ a finite subgroup of ${\rm SL}_2(k)$
such that $\#G$ is non-zero in $k$. Then $\CM S^G=\add S$ holds, and the Auslander algebra
$\End_{S^G}(S)$ is isomorphic to the \emph{skew group ring} $S*G$.
This is a free $S$-module with basis $G$, and the multiplication is given by $(sg)(s'g')=sg(s')gg'$
for $s,s'\in S$ and $g,g'\in G$. 
Thus the AR quiver of $\CM S^G$ coincides with the \emph{Gabriel quiver} of $S*G$,
and hence with the \emph{McKay quiver} of $G$,
which is the double of an extended Dynkin quiver. This is called \emph{algebraic McKay correspondence}.

On the other hand, the dual graph of the exceptional curves in the minimal resolution $X$ of
the singularity $\Spec S^G$ is a Dynkin graph. This is called \emph{geometric McKay correspondence}. 
There is a geometric construction of CM $S^G$-modules using $X$ \cite{AV}, 
which is a prototype of non-commutative crepant resolutions \cite{V1,V2}.
\end{example}

\section{Tilting and cluster tilting}

\subsection{Tilting theory}
Tilting theory is a Morita theory for triangulated categories.
It has an origin in Bernstein-Gelfand-Ponomarev reflection for quiver representations,
and established by works of Brenner-Butler, Happel-Ringel, Rickard, Keller and others (see e.g.\ \cite{AHK}).
The class of silting objects was introduced to complete the class of tilting objects
in the study of t-structures \cite{KV} and mutation \cite{AI}.

\begin{definition}
Let $\TT$ be a triangulated category.
A full subcategory of $\TT$ is \emph{thick} if it is closed under cones, $[\pm1]$ and direct summands.
We call an object $T\in\TT$ \emph{tilting} (resp.\ \emph{silting}) if $\Hom_{\TT}(T,T[i])=0$ holds
for all integers $i\neq0$ (resp.\ $i>0$), and the smallest thick subcategory of $\TT$
containing $T$ is $\TT$.
\end{definition}

The principal example of tilting objects appears in $\KKK^{\bo}(\proj\Lambda)$ for a ring $\Lambda$.
It has a tilting object given by the stalk complex $\Lambda$ concentrated in degree zero.
Conversely, any triangulated category with a tilting object is triangle equivalent to $\KKK^{\bo}(\proj\Lambda)$
under mild assumptions (see \cite{Ki2} for a detailed proof).

\begin{theorem}\cite{Ke1}\label{tilting theorem}
Let $\TT$ be an algebraic triangulated category and $T\in\TT$ a tilting object.
If $\TT$ is idempotent complete, then there is a triangle equivalence
$\TT\simeq\KKK^{\bo}(\proj\End_{\TT}(T))$ sending $T$ to $\End_{\TT}(T)$.
\end{theorem}

As an application, one can deduce Rickard's fundamental Theorem \cite{Ric1}, characterizing
when two rings are derived equivalent in terms of tilting objects.
Another application is the following converse of Proposition \ref{happel ZQ}(b).

\begin{example}\label{ZQ is D(kQ)}
Let $\TT$ be a $k$-linear Hom-finite Krull-Schmidt algebraic triangulated category over
an algebraically closed field $k$.
If the AR quiver of $\TT$ is $\Z Q$ for a Dynkin quiver $Q$, then $\TT$ has a tilting object
$T=\bigoplus_{i\in Q_0}(0,i)$ for $(0,i)\in\Z\times Q_0=(\Z Q)_0=\ind\TT$.
Thus there is a triangle equivalence $\TT\simeq\DDD^{\bo}(\mod kQ^{\op})$.
\end{example}

The following is the first main problem we will discuss in this paper.

\begin{problem}\label{first problem}
Find a $G$-graded Iwanaga-Gorenstein ring $\Lambda$ such that there is a triangle equivalence
$\underline{\CM}^G\Lambda\simeq\KKK^{\bo}(\proj\Gamma)$ for some ring $\Gamma$.
Equivalently (by Theorem \ref{tilting theorem}), find a $G$-graded Iwanaga-Gorenstein ring $\Lambda$ such that there is a tilting object in $\underline{\CM}^G\Lambda$.
\end{problem}

\subsection{Cluster tilting and higher Auslander-Reiten theory}\label{sec: Cluster tilting theory}
The notion of cluster tilting appeared naturally in a context of higher Auslander-Reiten theory \cite{I3}.
It also played a central role in categorification of cluster algebras \cite{FZ}
by using cluster categories, a new class of triangulated categories introduced in \cite{BMRRT}, and preprojective algebras \cite{GLS}.
Here we explain only the minimum necessary background for the aim of this paper.

Let $\Lambda$ be a finite dimensional $k$-algebra with $\gl\Lambda\le d$.
Then $\DDD^{\bo}(\mod\Lambda)$ has a Serre functor $\nu$
by Example \ref{serre functor of derived category}.
Using the \emph{higher AR translation} $\nu_d:=\nu\circ[-d]$ of $\DDD^{\bo}(\mod\Lambda)$,
the \emph{orbit category} $\CCC^{\circ}_d(\Lambda)=\DDD^{\bo}(\mod\Lambda)/\nu_d$ is defined.
It has the same objects as $\DDD^{\bo}(\mod\Lambda)$, and the morphism space is given by
\[\Hom_{\CCC^{\circ}_d(\Lambda)}(X,Y)=\bigoplus_{i\in\Z}\Hom_{\DDD^{\bo}(\mod\Lambda)}(X,\nu_d^i(Y)),\]
where the composition is defined naturally.
In general, $\CCC^{\circ}_d(\Lambda)$ does not have a natural structure of a triangulated category.
The \emph{$d$-cluster category} of $\Lambda$ is a triangulated category $\CCC_d(\Lambda)$
containing $\CCC^{\circ}_d(\Lambda)$ as a full subcategory such that the composition
$\DDD^{\bo}(\mod\Lambda)\to\CCC^{\circ}_d(\Lambda)\subset\CCC_d(\Lambda)$ is a triangle functor.
It was constructed in \cite{BMRRT} for hereditary case where $\CCC_d(\Lambda)=\CCC^{\circ}_d(\Lambda)$ holds, and in \cite{Ke2,Ke3,A,Gu} for general case by using a DG enhancement of $\DDD^{\bo}(\mod\Lambda)$.

We say that $\Lambda$ is \emph{$\nu_d$-finite} if $H^0(\nu_d^{-i}(\Lambda))=0$ holds for $i\gg0$.
This is automatic if $\gl\Lambda<d$.
In the hereditary case $d=1$, $\Lambda$ is $\nu_1$-finite if and only if it is representation-finite.
The following is a basic property of $d$-cluster categories.

\begin{theorem}\cite{A,Gu}
Let $\Lambda$ be a finite dimensional $k$-algebra with $\gl\Lambda\le d$.
Then $\Lambda$ is $\nu_d$-finite if and only if $\CCC_d(\Lambda)$ is Hom-finite.
In this case, $\CCC_d(\Lambda)$ is a $d$-Calabi-Yau triangulated category.
\end{theorem}

Thus, if $\Lambda$ is $\nu_d$-finite, then $\CCC_d(\Lambda)$ never has a tilting object.
But the object $\Lambda$ in $\CCC_d(\Lambda)$ still enjoys a similar property to tilting objects.
Now we recall the following notion, introduced in \cite{I1} as a
\emph{maximal $(d-1)$-orthogonal} subcategory. 

\begin{definition}\cite{I1}
Let $\TT$ be a triangulated or exact category and $d\ge1$. We call a full subcategory $\CC$ of $\TT$
\emph{$d$-cluster tilting} if $\CC$ is a functorially finite subcategory of $\TT$ such that
\begin{eqnarray*}
\CC&=&\{X\in\TT\mid\forall i\in\{1,2,\ldots,d-1\}\ \Ext_{\TT}^i(\CC,X)=0\}\\
&=&\{X\in\TT\mid\forall i\in\{1,2,\ldots,d-1\}\ \Ext_{\TT}^i(X,\CC)=0\}.
\end{eqnarray*}
We call an object $T\in\TT$ \emph{$d$-cluster tilting} if $\add T$ is a $d$-cluster tilting subcategory.
\end{definition}

If $\TT$ has a Serre functor $\SSS$, then it is easy to show $(\SSS\circ[-d])(\CC)=\CC$.
Thus it is natural in our setting $\TT=\DDD^{\bo}(\mod\Lambda)$ to consider the full subcategory
\begin{equation}\label{define Ud}
\UUU_d(\Lambda):=\add\{\nu_d^i(\Lambda)\mid i\in\Z\}\subset\DDD^{\bo}(\mod\Lambda).
\end{equation}
Equivalently, $\UUU_d(\Lambda)=\pi^{-1}(\add\pi\Lambda)$ for the functor
$\pi\colon \DDD^{\bo}(\mod\Lambda)\to\CCC_d(\Lambda)$.
In the hereditary case $d=1$, $\UUU_1(\Lambda)=\DDD^{\bo}(\mod\Lambda)$ holds if $\Lambda$ is 
representation-finite, and otherwise $\UUU_1(\Lambda)$ is the connected component of
the AR quiver of $\DDD^{\bo}(\mod\Lambda)$ containing $\Lambda$. 
This observation is generalized as follows.

\begin{theorem}\cite{A,I4}\label{C_d has d-CT}
Let $\Lambda$ be a finite dimensional $k$-algebra with $\gl\Lambda\le d$. If $\Lambda$ is $\nu_d$-finite,
then $\CCC_d(\Lambda)$ has a $d$-cluster tilting object $\Lambda$, and
$\DDD^{\bo}(\mod\Lambda)$ has a $d$-cluster tilting subcategory $\UUU_d(\Lambda)$.
\end{theorem}

We define a full subcategory of $\DDD^{\bo}(\mod\Lambda)$ by
\[\DDD^{d\Z}(\mod\Lambda)=\{X\in\DDD^{\bo}(\mod\Lambda)\mid\forall i\in\Z\setminus d\Z,\ H^i(X)=0\}.\]
If $\gl\Lambda\le d$, then any object in $\DDD^{d\Z}(\mod\Lambda)$ is isomorphic to a finite direct sum of
$X[di]$ for some $X\in\mod\Lambda$ and $i\in\Z$.
This generalizes Proposition \ref{happel ZQ}(a) for hereditary algebras,
and motivates the following definition.

\begin{definition}\cite{HIO}
Let $d\ge1$. A finite dimensional $k$-algebra $\Lambda$ is called \emph{$d$-hereditary}
if $\gl\Lambda\le d$ and $\UUU_d(\Lambda)\subset\DDD^{d\Z}(\mod\Lambda)$.
\end{definition}

It is clear that $1$-hereditary algebras are precisely hereditary algebras.
We have the following dichotomy of $d$-hereditary algebras.

\begin{theorem}\cite{HIO}
Let $\Lambda$ be a ring-indecomposable finite dimensional $k$-algebra with $\gl\Lambda\le d$.
Then $\Lambda$ is $d$-hereditary if and only if either (i) or (ii) holds:
\begin{enumerate}[\rm(i)]
\item There exists a $d$-cluster tilting object in $\mod\Lambda$.
\item $\nu_d^{-i}(\Lambda)\in\mod\Lambda$ holds for any $i\ge0$.
\end{enumerate}
\end{theorem}

When $d=1$, the above (i) holds if and only if $\Lambda$ is representation-finite,
and the above (ii) holds if and only if $\Lambda$ is $d$-representation-infinite.

\begin{definition}\label{define dichotomy}
Let $\Lambda$ be a $d$-hereditary algebra.
We call $\Lambda$ \emph{$d$-representation-finite} if the above (i) holds, and
\emph{$d$-representation-infinite} if the above (ii) holds.
\end{definition}

\begin{example}\label{example of d-hereditary}
\begin{enumerate}[\rm(a)]
\item Let $\Lambda=kQ$ for a connected acyclic quiver $Q$. Then $\Lambda$ is $1$-representation-finite if $Q$ is Dynkin, 
and $1$-representation-infinite otherwise.
\item Let $X$ be a smooth projective variety with $\dim X=d$, and $T\in\coh X$ a tilting object in
$\DDD^{\bo}(\coh X)$. Then $\Lambda=\End_X(T)$ always satisfies $\gl\Lambda\ge d$.
If the equality holds, then $\Lambda$ is $d$-representation-infinite \cite{BuH}.
\item There is a class of finite dimensional $k$-algebras called \emph{Fano algebras} \cite{M,MM}
in non-commutative algebraic geometry. So-called \emph{extremely Fano algebras} $\Lambda$
with $\gl\Lambda=d$ are $d$-representation-infinite.
\end{enumerate}
\end{example}

It is known that $d$-cluster tilting subcategories of a triangulated (resp. exact) category $\TT$
enjoy various properties which should be regarded as higher analogues of those of $\TT$.
For example, they have \emph{almost split $(d+2)$-angles} by \cite{IYo}
(resp. \emph{$d$-almost split sequences} by \cite{I1}), and structures of
\emph{$(d+2)$-angulated categories} by \cite{GKO} (resp. \emph{$d$-abelian categories} by \cite{Jas}).
These motivate the following definition.

\begin{definition}[cf.\ Definition \ref{define CM-finite}]\label{define d-CM-finite}
An Iwanaga-Gorenstein ring $\Lambda$ is called \emph{$d$-CM-finite} if 
there exists a $d$-cluster tilting object $M$ in $\CM\Lambda$.
In this case, we call $\End_\Lambda(M)$ the \emph{$d$-Auslander algebra} and
$\underline{\End}_\Lambda(M)$ the \emph{stable $d$-Auslander algebra}.
\end{definition}

$1$-CM-finiteness coincides with classical CM-finiteness since
$1$-cluster tilting objects are precisely additive generators. 
\emph{$d$-Auslander correspondence} gives a characterization of a certain nice class of
algebras with finite global dimension as $d$-Auslander algebras \cite{I2}.
As a special case, it gives a connection with non-commutative crepant resolutions (NCCRs) of
Van den Bergh \cite{V2}.
Recall that a reflexive module $M$ over a Gorenstein ring $R$ gives an \emph{NCCR}
$\End_R(M)$ of $R$ if $\End_R(M)$ is a non-singular $R$-order
(see Definition \ref{example of IG}(c)).

\begin{theorem}\cite{I2}
Let $R$ be a Gorenstein ring with $\dim R=d+1$. Assume $M\in\CM R$ has $R$ as a direct summand.
Then $M$ is a $d$-cluster tilting object in $\CM R$ if and only if $M$ gives an NCCR of $R$ and $R$ satisfies Serre's (R$_d$) condition.
\end{theorem}

The following generalizes Example \ref{quotient singularity2}.

\begin{example}\cite{I1,V2}\label{quotient singularity}
Let $S=k[[x_0,\ldots,x_d]]$ be a formal power series ring and $G$ a finite subgroup of ${\rm SL}_{d+1}(k)$
such that $\# G$ is non-zero in $k$. Then the $S^G$-module $S$ gives an NCCR $\End_{S^G}(S)=S*G$
of $S^G$. If $S^G$ has at worst an isolated singularity, then $S$ is
a $d$-cluster tilting object in $\CM S^G$, and hence $S^G$ is $d$-CM-finite
with the $d$-Auslander algebra $S*G$.
As in Example \ref{quotient singularity2}, the quiver of $\add S$ coincides with
the Gabriel quiver of $S*G$ and with the McKay quiver of $G$.
\end{example}

The following is the second main problem we will discuss in this paper.

\begin{problem}\label{second problem}
Find a $d$-CM-finite Iwanaga-Gorenstein ring.
More strongly (by Theorem \ref{C_d has d-CT}), find an Iwanaga-Gorenstein ring $\Lambda$
such that there is a triangle equivalence
$\underline{\CM}\Lambda\simeq\CCC_d(\Gamma)$ for some algebra $\Gamma$.
\end{problem}

We refer to \cite{EH,Ber} for some necessary conditions for $d$-CM-finiteness.
Besides results in this paper, a number of examples of NCCRs have been found, see e.g.\ \cite{L,W,SV} and references therein.

It is natural to ask how the notion of $d$-CM-finiteness is related to CM-tameness
(e.g.\ \cite{BD}) and also the representation type of homogeneous coordinate rings
of projective varieties (e.g.\ \cite{FM}).

\section{Results in dimension 0 and 1}

\subsection{Dimension zero}

In this subsection, we consider finite dimensional Iwanaga-Gorenstein algebras.
We start with a classical result due to Happel \cite{H}.
Let $\Lambda$ be a finite dimensional $k$-algebra.
The \emph{trivial extension algebra} of $\Lambda$ is $T(\Lambda)=\Lambda\oplus D\Lambda$,
where the multiplication is given by $(\lambda,f)(\lambda',f')=(\lambda\lambda',\lambda f'+f\lambda')$
for $(\lambda,f),(\lambda',f')\in T(\Lambda)$. This is clearly a self-injective $k$-algebra, and has a
$\Z$-grading given by $T(\Lambda)_0=\Lambda$, $T(\Lambda)_1=D\Lambda$ and
$T(\Lambda)_i=0$ for $i\neq0,1$.

\begin{theorem}\cite{H}\label{happel repetitive}
Let $\Lambda$ be a finite dimensional $k$-algebra with $\gl\Lambda<\infty$.
Then $\underline{\mod}^{\Z}T(\Lambda)$ has a tilting object $\Lambda$ such that $\underline{\End}^{\Z}_{T(\Lambda)}(\Lambda)\simeq \Lambda$, and there is a triangle equivalence
\begin{equation}\label{happel repetitive2}
\underline{\mod}^{\Z}T(\Lambda)\simeq\DDD^{\bo}(\mod \Lambda).
\end{equation}
\end{theorem}

As an application, it follows from Gabriel's Theorem \ref{gabriel} and covering theory
that $T(kQ)$ is representation-finite for any Dynkin quiver $Q$.
More generally, a large family of representation-finite self-injective algebras was constructed
from Theorem \ref{happel repetitive}. See a survey article \cite{S}.

Recently, Theorem \ref{happel repetitive} was generalized to a large class of $\Z$-graded
self-injective algebras $\Lambda$.
For $X\in\mod^{\Z}\Lambda$, let $X_{\ge0}=\bigoplus_{i\ge0}X_i$.

\begin{theorem}\label{dimension 0}\cite{Ya}
Let $\Lambda=\bigoplus_{i\ge0}\Lambda_i$ be a $\Z$-graded finite dimensional self-injective
$k$-algebra such that $\gl\Lambda_0<\infty$.
Then $\underline{\mod}^{\Z}\Lambda$ has a tilting object $T=\bigoplus_{i>0}\Lambda(i)_{\ge0}$, and
there is a triangle equivalence $\underline{\mod}^{\Z}\Lambda\simeq\KKK^{\bo}(\proj\underline{\End}^{\Z}_\Lambda(T))$.
\end{theorem}

If  $\soc\Lambda\subset\Lambda_a$ for some $a\in\Z$, then
$\underline{\End}^{\Z}_\Lambda(T)$ has a simple description
{\small\[\underline{\End}^{\Z}_\Lambda(T)\simeq
\begin{bmatrix}
\Lambda_0&0&\cdots&0&0\\
\Lambda_1&\Lambda_0&\cdots&0&0\\
\vdots&\vdots&\ddots&\vdots&\vdots\\
\Lambda_{a-2}&\Lambda_{a-3}&\cdots&\Lambda_0&0\\
\Lambda_{a-1}&\Lambda_{a-2}&\cdots&\Lambda_1&\Lambda_0\end{bmatrix}.\]}
For example, if $\Lambda=k[x]/(x^{a+1})$ with $\deg x=1$, then $\underline{\End}^{\Z}_\Lambda(T)$
is the path algebra $k\A_a$ of the quiver of type $A_{a}$.

We end this subsection with posing the following open problem.

\begin{problem}\label{tilting for IG}
Let $\Lambda=\bigoplus_{i\ge0}\Lambda_i$ be a $\Z$-graded finite dimensional Iwanaga-Gorenstein algebra.
When does $\underline{\CM}^{\Z}\Lambda$ have a tilting object?
\end{problem}

Recently, it was shown in \cite{LZ} and \cite{KMY} independently that 
if $\Lambda=\bigoplus_{i\ge0}\Lambda_i$ is a $\Z$-graded finite dimensional 1-Iwanaga-Gorenstein
algebra satisfying $\gl\Lambda_0<\infty$, then the stable category $\underline{\CM}^{\Z}\Lambda$
has a silting object.
We will see some other results in Section \ref{section: Truncated preprojective algebras}.
We refer to \cite{DI} for some results on Problem \ref{second problem}.

\subsection{Dimension one}
In this subsection, we consider a $\Z$-graded Gorenstein ring $R=\bigoplus_{i\ge0}R_i$
with $\dim R=1$ such that $R_0$ is a field.
Let $S$ be the multiplicative set of all homogeneous non-zerodivisors of $R$, and $K=RS^{-1}$ the
$\Z$-graded total quotient ring. Then there exists a positive integer $p$ such that $K(p)\simeq K$ as
$\Z$-graded $R$-modules.
In this setting, we have the following result (see Definitions \ref{define CM_0} and \ref{define a-invariant}
for $\CM_0^{\Z}R$ and the $a$-invariant).

\begin{theorem}\cite{BIY}\label{biy}
Let $R=\bigoplus_{i\ge0}R_i$ be a $\Z$-graded Gorenstein ring with $\dim R=1$ such that $R_0$ is a field,
and $a$ the $a$-invariant of $R$.
\begin{enumerate}[\rm(a)]
\item Assume $a\ge0$. Then $\underline{\CM}_0^{\Z}R$ has a tilting object
$T=\bigoplus_{i=1}^{a+p}R(i)_{\ge0}$, and there is a triangle equivalence
$\underline{\CM}_0^{\Z}R\simeq\KKK^{\bo}(\proj\underline{\End}^{\Z}_R(T))$.
\item Assume $a<0$. Then $\underline{\CM}_0^{\Z}R$ has a silting object $\bigoplus_{i=1}^{a+p}R(i)_{\ge0}$.
Moreover, it has a tilting object if and only if $R$ is regular.
\end{enumerate}
\end{theorem}

An important tool in the proof is Theorem \ref{orlov embedding}.
The endomorphism algebra of $T$ above has the following description.
{\small\begin{equation*}
\underline{\End}^{\Z}_R(T)=
\begin{bmatrix}
R_0&0&\cdots&0&0&
0&0&\cdots&0&0\\
R_1&R_0&\cdots&0&0&
0&0&\cdots&0&0\\
\vdots&\vdots&\ddots&\vdots&\vdots
&\vdots&\vdots&\cdots&\vdots&\vdots\\
R_{a-2}&R_{a-3}&\cdots&R_0&0&
0&0&\cdots&0&0\\
R_{a-1}&R_{a-2}&\cdots&R_1&R_0&
0&0&\cdots&0&0\\
K_a&K_{a-1}&\cdots&K_2&K_1&
K_0&K_{-1}&\cdots&K_{2-p}&K_{1-p}\\
K_{a+1}&K_a&\cdots&K_3&K_2&
K_{1}&K_0&\cdots&K_{3-p}&K_{2-p}\\
\vdots&\vdots&\vdots&\vdots&\vdots
&\vdots&\vdots&\ddots&\vdots&\vdots\\
K_{a+p-2}&K_{a+p-3}&\cdots&K_p&K_{p-1}&
K_{p-2}&K_{p-3}&\cdots&K_{0}&K_{-1}\\
K_{a+p-1}&K_{a+p-2}&\cdots&K_{p+1}&K_p&
K_{p-1}&K_{p-2}&\cdots&K_{1}&K_0
\end{bmatrix}.\end{equation*}}
As an application, we obtain the following graded version of Proposition \ref{AR quiver of simple singularity}(b).

\begin{example}\label{simple curve singularity}
Let $R=k[x,y]/(f)$ be a simple singularity \eqref{simple equation} with $\dim R=1$ and the
grading given by the list below. Then there is a triangle equivalence
$\underline{\CM}^{\Z}R\simeq\DDD^{\bo}(\mod kQ)$, 
where $Q$ is the Dynkin quiver in the list below.
In particular, the AR quiver of $\underline{\CM}^{\Z}R$ is $\Z Q^{\op}$ \cite{Ar}.
{\small\[\begin{array}{|c||c|c|c|c|c|c|c|}
\hline
R&A_{2n-1}&A_{2n}&D_{2n}&D_{2n+1}&E_6&E_7&E_8\\ \hline\hline
(\deg x,\deg y)&(1,n)&(2,2n+1)&(1,n-1)&(2,2n-1)&(3,4)&(2,3)&(3,5)\\ \hline
Q&D_{n+1}&A_{2n}&D_{2n}&A_{4n-1}&E_6&E_7&E_8\\ \hline
\end{array}\]}
This gives a conceptual proof of the classical result that simple singularities in dimension 1
are CM-finite \cite{J,DR,GK}.
\end{example}

In the special case below, one can construct a different tilting object, whose endomorphism
algebra is 2-representation-finite (Definition \ref{define dichotomy}).
This is closely related to the 2-cluster tilting object constructed in \cite{BIKR}.

\begin{theorem}\cite{HI2}
Let $R=k[x,y]/(f)$ be a hypersurface singularity with $f=f_1f_2\cdots f_n$ for linear forms $f_i$
and $\deg x=\deg y=1$. Assume that $R$ is reduced.
\begin{enumerate}[\rm(a)]
\item $\underline{\CM}^{\Z}R$ has a tilting object $U=\bigoplus_{i=1}^n(k[x,y]/(f_1f_2\cdots f_i)\oplus k[x,y]/(f_1f_2\cdots f_i)(1))$.
\item $\underline{\End}_R^{\Z}(U)$ is a 2-representation-finite algebra. It is the Jacobian algebra of a certain quiver with potential.
\end{enumerate}
\end{theorem}

We refer to \cite{DL,JKS,Ge} for other results in dimension one.

\section{Preprojective algebras}

\subsection{Classical preprojective algebras}
Preprojective algebras are widely studied objects with various applications,
e.g.\ cluster algebras \cite{GLS}, quantum groups \cite{KS,Lu}, quiver varieties \cite{Na}.
Here we discuss a connection to CM representations.

Let $Q$ be an acyclic quiver, and $\overline{Q}$ the double of $Q$ obtained by
adding an inverse arrow $a^*\colon j\to i$ for each arrow $a\colon i\to j$ in $Q$.
The \emph{preprojective algebra} of $Q$ is the factor algebra of the path algebra $k\overline{Q}$ defined by
\begin{equation}\label{define pi}
\Pi=k\overline{Q}/(\sum_{a\in Q_1}(aa^*-a^*a)).
\end{equation}
We regard $\Pi$ as a $\Z$-graded algebra by $\deg a=0$ and $\deg a^*=1$ for any $a\in Q_1$.
Clearly $\Pi_0=kQ$ holds. Moreover $\Pi_1=\Ext^1_{kQ}(D(kQ),kQ)$ as a $kQ$-bimodule,
and $\Pi$ is isomorphic to the tensor algebra $T_{kQ}\Ext^1_{kQ}(D(kQ),kQ)$.
Thus the $kQ$-module $\Pi_i$ is isomorphic to the \emph{preprojective $kQ$-module} $H^0(\tau^{-i}(kQ))$,
where $\tau=\nu\circ[-1]$ is the AR translation.
This is the reason why $\Pi$ is called the preprojective algebra. Moreover, for the category $\UUU_1(kQ)$ defined in \eqref{define Ud}, there is an equivalence
\begin{equation}\label{U1 is proj}
\UUU_1(kQ)=\add\{\tau^{-i}(kQ)\mid i\in\Z\}\simeq\proj^{\Z}\Pi
\end{equation}
given by $X\mapsto\bigoplus_{i\in\Z}\Hom_{\UUU_1(kQ)}(kQ,\tau^{-i}(X))$,
which gives the following trichotomy.
\[\begin{array}{|c||c|c|c|}\hline
Q&\mbox{Dynkin}&\mbox{extended Dynkin}&\mbox{else}\\ \hline\hline
kQ&\mbox{representation-finite}&\mbox{representation-tame}&\mbox{representation-wild}\\ \hline
\dim_k\Pi_i&\dim_k\Pi<\infty&\mbox{linear growth}&\mbox{exponential growth}\\ \hline
\end{array}\]
It was known in 1980s that, $\Pi$ in the extended Dynkin case has
a close connection to simple singularities.

\begin{theorem}\cite{Au1,GL1,GL2,RV1}\label{Pi for extended dynkin}
Let $\Pi$ be a preprojective algebra of an extended Dynkin quiver $Q$, $e$ the vertex $\circ$ in
\eqref{dynkin}, and $R=e\Pi e$.
\begin{enumerate}[\rm(a)]
\item $R$ is the simple singularity $k[x,y,z]/(f)$ in dimension 2 with the induced $\Z$-grading below, where $p$ in type $A_n$ is the number of clockwise arrows in $Q$.
{\small\[\begin{array}{|c||c|c|c|}
\hline
Q,R&A_{n}&D_{2n}&D_{2n+1}\\ \hline\hline
f&x^{n+1}-yz&x(y^2+x^{n-1}y)+z^2&x(y^2+x^{n-1}z)+z^2\\ \hline
(\deg x,\deg y,\deg z)&(1,p,n+1-p)&(2,2n-2,2n-1)&(2,2n-1,2n)\\ \hline\hline
Q,R&E_6&E_7&E_8\\ \hline\hline
f&x^2z+y^3+z^2&x^3y+y^3+z^2&x^5+y^3+z^2\\ \hline
(\deg x,\deg y,\deg z)&(3,4,6)&(4,6,9)&(6,10,15)\\ \hline
\end{array}
\]}
\item $\Pi e$ is an additive generator of $\CM R$ and satisfies $\End_R(\Pi e)=\Pi$.
Therefore $R$ is CM-finite with an Auslander algebra $\Pi$.
\item $\Pi$ is Morita equivalent to the skew group ring $k[u,v]*G$ for a finite subgroup $G$ of $\SL_2(k)$
if $k$ is sufficiently large (cf.\ Example \ref{quotient singularity2}).
\end{enumerate}
\end{theorem}

By (b) and \eqref{U1 is proj} above, there are equivalences $\CM^{\Z}R\simeq\proj^{\Z}\Pi\simeq\UUU_1(kQ)\subset\DDD^{\bo}(\mod kQ)$. 
Thus the AR quivers of $\CM^{\Z}R$ and $\underline{\CM}^{\Z}R$
are given by $\Z Q^{\op}$ and $\Z(Q^{\op}\backslash\{e\})$ respectively.
Now the following result follows from Example \ref{ZQ is D(kQ)}.

\begin{corollary}\label{simple surface singularity}
Under the setting in Theorem \ref{Pi for extended dynkin}, there is a triangle equivalence
$\underline{\CM}^{\Z}R\simeq\DDD^{\bo}(\mod kQ/(e))$.
\end{corollary}

Two other proofs were given in \cite{KST1},
one uses explicit calculations of $\Z$-graded matrix factorizations,
and the other uses Theorem \ref{orlov embedding}.
In Theorem \ref{air theorem} below, we deduce Corollary \ref{simple surface singularity} from
a general result on higher preprojective algebras.
We refer to \cite{KST2,LP} for results for some other hypersurfaces in dimension 2.

\subsection{Higher preprojective algebras}\label{section: d-representation infinite case}

There is a natural analogue of preprojective algebras for finite dimensional algebras with finite global dimension. 

\begin{definition}\cite{IO}
Let $\Lambda$ be a finite dimensional $k$-algebra with $\gl\Lambda\le d$.
We regard the highest extension $\Ext^d_\Lambda(D\Lambda,\Lambda)$ as a $\Lambda$-bimodule naturally, 
and define the \emph{$(d+1)$-preprojective algebra} as the tensor algebra
\[\Pi_{d+1}(\Lambda)=T_\Lambda\Ext^d_\Lambda(D\Lambda,\Lambda).\]
\end{definition}

This is the 0-th cohomology of the \emph{Calabi-Yau completion} of $\Lambda$ \cite{Ke3}.
For example, for an acyclic quiver $Q$, $\Pi_2(kQ)$ is the preprojective algebra \eqref{define pi}.

The algebra $\Pi=\Pi_{d+1}(\Lambda)$ has an alternative description in terms of the higher AR translation
$\nu_d=\nu\circ[-d]$ of $\DDD^{\bo}(\mod\Lambda)$. The $\Z$-grading on $\Pi$ is given by
\[\Pi_i=\Ext^d_\Lambda(D\Lambda,\Lambda)^{\otimes_\Lambda i}=
\Hom_{\DDD^{\bo}(\mod\Lambda)}(\Lambda,\nu_d^{-i}(\Lambda))\]
for $i\ge0$. Thus there is an isomorphism $\Pi\simeq\End_{\CCC_d(\Lambda)}(\Lambda)$ and an equivalence
\begin{equation}\label{U is proj Pi}
\UUU_d(\Lambda)\simeq\proj^{\Z}\Pi
\end{equation}
given by $X\mapsto\bigoplus_{i\in\Z}\Hom_{\UUU_d(\Lambda)}(\Lambda,\nu_d^{-i}(X))$.
In particular, $\Pi$ is finite dimensional if and only if $\Lambda$ is $\nu_d$-finite.

We see below that $\Pi_{d+1}(\Lambda)$ enjoys nice homological properties if $\Lambda$ is $d$-hereditary.

\begin{definition}[cf.\ \cite{Gi}]
Let $\Gamma=\bigoplus_{i\ge0}\Gamma_i$ be a $\Z$-graded $k$-algebra.
We denote by $\Gamma^{\rm e}=\Gamma^{\op}\otimes_k\Gamma$ the enveloping algebra of $\Gamma$.
We say that $\Gamma$ is a \emph{$d$-Calabi-Yau algebra of $a$-invariant $a$} (or \emph{Gorenstein parameter $-a$}) if $\Gamma$ belongs to $\KKK^{\bo}(\proj^{\Z}\Gamma^{\rm e})$
and $\RHom_{\Gamma^{\rm e}}(\Gamma,\Gamma^{\rm e})(a)[d]\simeq\Gamma$ holds in
$\DDD(\Mod^{\Z}\Gamma^{\rm e})$.
\end{definition}

For example, the $\Z$-graded polynomial algebra $k[x_1,\ldots,x_{d}]$ with $\deg x_i=a_i$ is a $d$-Calabi-Yau algebra of $a$-invariant $-\sum_{i=1}^na_i$.

Now we give a homological characterization of the $(d+1)$-preprojective algebras of
$d$-representation-infinite algebras (Definition \ref{define dichotomy}) as the explicit correspondence.

\begin{theorem}\cite{Ke3,MM,AIR}\label{MM correspondence}
There exists a bijection between the set of isomorphism classes of $d$-representation-infinite algebras $\Lambda$ and the set of isomorphism classes of $(d+1)$-Calabi-Yau algebras $\Gamma$ of
$a$-invariant $-1$. It is given by $\Lambda\mapsto\Pi_{d+1}(\Lambda)$ and $\Gamma\mapsto\Gamma_0$.
\end{theorem}

Note that $\Gamma$ above is usually non-noetherian. If $\Gamma$ is right graded coherent, then for the category $\qgr\Gamma$ defined in \eqref{define qgr}, there is a triangle equivalence \cite{M}
\begin{equation}\label{minamoto equivalence}
\DDD^{\bo}(\mod\Lambda)\simeq\DDD^{\bo}(\qgr\Gamma).
\end{equation}
Applying Theorem \ref{MM correspondence} for $d=1,2$, we obtain the following observations
(see \cite{V3} for a structure theorem of (ungraded) Calabi-Yau algebras).

\begin{example}
Let $k$ be an algebraically closed field.
\begin{enumerate}[\rm(a)]
\item (cf.\ \cite{Bo1}) $2$-Calabi-Yau algebras of $a$-invariant $-1$ are precisely the preprojective
algebras of disjoint unions of non-Dynkin quivers.
\item (cf.\ \cite{Bo1,HI1}) $3$-Calabi-Yau algebras of $a$-invariant $-1$ are precisely the Jacobian algebras of
quivers with `good' potential with cuts.
\end{enumerate}
\end{example}

The setting of our main result is the following.

\begin{assumption}\label{air assumption}
Let $\Gamma$ be a $(d+1)$-Calabi-Yau algebras of $a$-invariant $-1$.
Equivalently by Theorem \ref{MM correspondence}, $\Gamma$ is a $(d+1)$-preprojective algebra
of some $d$-representation-infinite algebra.
We assume that the following conditions hold for $\Lambda=\Gamma_0$.
\begin{enumerate}[\rm(i)]
\item $\Gamma$ is a noetherian ring, $e\in\Lambda$ is an idempotent and $\dim_k(\Gamma/(e))<\infty$. 
\item $e\Lambda(1-e)=0$.
\end{enumerate}
\end{assumption}

For example, let $Q$ be an extended Dynkin quiver. If the vertex $\circ$ in \eqref{dynkin} is a sink,
then $\Gamma=\Pi_2(kQ)$ and $e=\circ$ satisfy Assumption \ref{air assumption}
by Theorem \ref{Pi for extended dynkin}.

Under Assumption \ref{air assumption}(i), let $R=e\Gamma e$. 
Then $R$ is a $(d+1)$-Iwanaga-Gorenstein ring, and the $(\Gamma, R)$-bimodule $\Gamma e$ plays an important role.
It is a CM $R$-module, and gives a $d$-cluster tilting object in $\CM R$. 
Moreover the natural morphism $\Gamma\to\End_R(\Gamma e)$ is an isomorphism.
Thus $R$ is $d$-CM-finite and has a $d$-Auslander algebra $\Gamma$.
The proof of these statements is parallel to Example \ref{quotient singularity}.

Regarding $\Gamma e$ as a $\Z$-graded $R$-module, we consider the composition
\[F\colon \DDD^{\bo}(\mod\Lambda/(e))\to\DDD^{\bo}(\mod\Lambda)\xrightarrow{-\Lotimes_{\Lambda}\Gamma e}\DDD^{\bo}(\mod^{\Z}R)\to\underline{\CM}^{\Z}R,\]
where the first functor is induced from the surjective morphism $\Lambda\to\Lambda/(e)$,
and the last functor is given by Theorem \ref{stable is singularity}.
Under Assumption \ref{air assumption}(ii), $F$ is shown to be a triangle equivalence.
A crucial step is to show that $F$ restricts to an equivalence
$\UUU_d(\Lambda/(e))\to\add\{\Gamma e(i)\mid i\in\Z\}$, which are $d$-cluster tilting subcategories of
$\DDD^{\bo}(\mod\Lambda/(e))$ and $\underline{\CM}^{\Z}R$ respectively (Theorem \ref{C_d has d-CT}). 
Similarly,  we obtain a triangle equivalence $\CCC_d(\Lambda/(e))\simeq\underline{\CM}R$
by using universality of $d$-cluster categories \cite{Ke2}.
As a summary, we obtain the following results.

\begin{theorem}\cite{AIR}\label{air theorem}
Under Assumption \ref{air assumption}(i), let $R=e\Gamma e$ and $\Lambda=\Gamma_0$.
\begin{enumerate}[\rm(a)]
\item $R$ is a $(d+1)$-Iwanaga-Gorenstein algebra, and $\Gamma e$ is a CM $R$-module.
\item $\Gamma e$ is a $d$-cluster tilting object in $\CM R$ and satisfies $\End_{R}(\Gamma e)=\Gamma$.
Thus $R$ is $d$-CM-finite and has a $d$-Auslander algebra $\Gamma$
\item If Assumption \ref{air assumption}(ii) is satisfied, then there exist triangle equivalences
\[\DDD^{\bo}(\mod\Lambda/(e))\simeq\underline{\CM}^{\Z}R\ \mbox{ and }\ \CCC_d(\Lambda/(e))\simeq\underline{\CM}R.\]
\end{enumerate}
\end{theorem}

Similar triangle equivalences were given in \cite{TV,KY} using different methods.
There is a connection between (c) and \eqref{minamoto equivalence} above via
Theorem \ref{orlov embedding}, see \cite{Am2}.

In the case $d=1$, the above (c) recovers Corollary \ref{simple surface singularity}
and a triangle equivalence $\underline{\CM}R\simeq\CCC_1(kQ/(e))$, which implies algebraic McKay correspondence in Example \ref{quotient singularity2}.
Motivated by Example \ref{quotient singularity} and Theorem \ref{Pi for extended dynkin}(c),
we consider the following.

\begin{example}\cite{AIR,U1}
Let $S=k[x_0,\ldots,x_d]$ be a polynomial algebra, and $G$ a finite subgroup of ${\rm SL}_{d+1}(k)$.
Then the skew group ring $\Gamma=S*G$ is a (ungraded) $(d+1)$-Calabi-Yau algebra.
Assume that $G$ is generated by the diagonal matrix ${\rm diag}(\zeta^{a_0},\ldots,\zeta^{a_d})$,
where $\zeta$ is a primitive $n$-th root of unity and $0\le a_j\le n-1$ for each $j$.
Then $\Gamma$ is presented by the McKay quiver of $G$, which has vertices $\Z/n\Z$,
and arrows $x_j\colon i\to i+a_j$ for each $i,j$. Define a $\Z$-grading on $\Gamma$ by
$\deg(x_j\colon i\to i+a_j)=0$ if $i<i+a_j$ as integers in $\{1,\ldots,n\}$, and $1$ otherwise.
Then $\Gamma$ is a $(d+1)$-Calabi-Yau algebra of $a$-invariant $-\sum_{0\le j\le d}a_j/n$.
Assume that this is $-1$, and let $e=e_n$.
Then Assumption \ref{air assumption} is satisfied, and $e\Gamma e=S^G$ holds.
Thus Theorem \ref{air theorem} gives triangle equivalences 
\[\DDD^{\bo}(\mod\Lambda/(e))\simeq\underline{\CM}^{\Z}S^G\ \mbox{ and }\ \CCC_d(\Lambda/(e))\simeq\underline{\CM}S^G.\]
Below we draw quivers for two cases (i) $n=d+1$ and $a_0=\cdots=a_d=1$, and
(ii) $d=2$, $n=5$ and $(a_0,a_1,a_2)=(1,2,2)$.
{\scriptsize
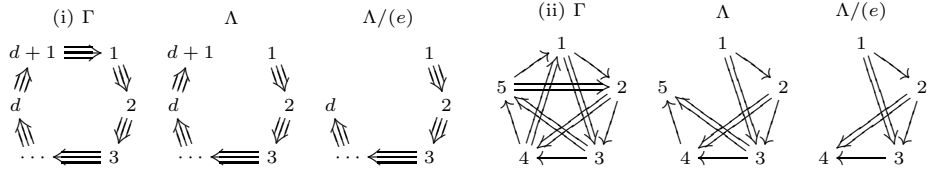
\[\begin{xy}
0;<2.2pt,0pt>:<0pt,2.2pt>:: 
(0,18) *+{{\rm(i)}\ \Gamma},
(7,12) *+{1} ="1",
(10,3) *+{2} ="2",
(7,-6) *+{3} ="3",
(-7,-6) *+{\cdots} ="4",
(-10,3) *+{d} ="5",
(-7,12) *+{d+1} ="6",
\ar@3{->} "1";"2"
\ar@3{->} "2";"3"
\ar@3{->} "3";"4"
\ar@3{->} "4";"5"
\ar@3{->} "5";"6"
\ar@3{->} "6";"1"
\end{xy}
\ \ \begin{xy}
0;<2.2pt,0pt>:<0pt,2.2pt>:: 
(0,18) *+{\Lambda},
(7,12) *+{1} ="1",
(10,3) *+{2} ="2",
(7,-6) *+{3} ="3",
(-7,-6) *+{\cdots} ="4",
(-10,3) *+{d} ="5",
(-7,12) *+{d+1} ="6",
\ar@3{->} "1";"2"
\ar@3{->} "2";"3"
\ar@3{->} "3";"4"
\ar@3{->} "4";"5"
\ar@3{->} "5";"6"
\end{xy}
\ \ \begin{xy}
0;<2.2pt,0pt>:<0pt,2.2pt>:: 
(0,18) *+{\Lambda/(e)},
(7,12) *+{1} ="1",
(10,3) *+{2} ="2",
(7,-6) *+{3} ="3",
(-7,-6) *+{\cdots} ="4",
(-10,3) *+{d} ="5",
\ar@3{->} "1";"2"
\ar@3{->} "2";"3"
\ar@3{->} "3";"4"
\ar@3{->} "4";"5"
\end{xy}\ \ \ \ 
\begin{xy}
0;<2pt,0pt>:<0pt,2pt>:: 
(0,21) *+{{\rm (ii)}\ \Gamma},
(0,15) *+{1} ="1",
(11.4,6.71) *+{2} ="2",
(7.1,-6.7) *+{3} ="3",
(-7.1,-6.7) *+{4} ="4",
(-11.4,6.71) *+{5} ="5",
\ar "1";"2"
\ar "2";"3"
\ar "3";"4"
\ar "4";"5"
\ar "5";"1"
\ar@{=>} "1";"3"
\ar@{=>} "2";"4"
\ar@{=>} "3";"5"
\ar@{=>} "4";"1"
\ar@{=>} "5";"2"
\end{xy}
\ \ \begin{xy}
0;<2pt,0pt>:<0pt,2pt>:: 
(0,21) *+{\Lambda},
(0,15) *+{1} ="1",
(11.4,6.71) *+{2} ="2",
(7.1,-6.7) *+{3} ="3",
(-7.1,-6.7) *+{4} ="4",
(-11.4,6.71) *+{5} ="5",
\ar "1";"2"
\ar "2";"3"
\ar "3";"4"
\ar "4";"5"
\ar@{=>} "1";"3"
\ar@{=>} "2";"4"
\ar@{=>} "3";"5"
\end{xy}
\ \ \begin{xy}
0;<2pt,0pt>:<0pt,2pt>:: 
(0,21) *+{\Lambda/(e)},
(0,15) *+{1} ="1",
(11.4,6.71) *+{2} ="2",
(7.1,-6.7) *+{3} ="3",
(-7.1,-6.7) *+{4} ="4",
\ar "1";"2"
\ar "2";"3"
\ar "3";"4"
\ar@{=>} "1";"3"
\ar@{=>} "2";"4"
\end{xy}\]
}
In (i), $\xymatrix@C1em{\ar@3{->}[r]&}$ shows $d+1$ arrows, $S^G$ is the Veronese subring
$S^{(d+1)}$ and $\Lambda$ is the Beilinson algebra. For $d=2$, we recover the triangle equivalence
$\CCC_2(kQ)\simeq\underline{\CM}S^G$ for
$Q=[\xymatrix@C1em{1\ar@<3pt>[r]\ar[r]\ar@<-3pt>[r]&2}]$ given in \cite{KR,KMV}.
\end{example}

Note that similar triangle equivalences are given in \cite{IT,U2,MU} for the skew group rings $S*G$ 
whose $a$-invariants are not equal to $-1$.

\begin{example}[Dimer models]
Let $G$ be a bipartite graph on a torus, and $G_0$ (resp.\ $G_1$, $G_2$)
the set of vertices (resp.\ edges, faces) of $G$. We associate a quiver with potential $(Q,W)$:
The underlying graph of $Q$ is the dual of the graph $G$, and
faces of $Q$ dual to white (resp.\ black) vertices are oriented clockwise
(resp.\ anti-clockwise).
Hence any vertex $v \in G_0$ corresponds to a cycle $c_v$ of $Q$.
Let $W=\sum_{v: \mbox{{\scriptsize white}}}c_v-\sum_{v: \mbox{{\scriptsize black}}}c_v$,
and $\Gamma$ the Jacobian algebra of $(Q,W)$.

Under the assumption that $G$ is \emph{consistent}, $\Gamma$ is a (ungraded) 3-Calabi-Yau algebra,
and for any vertex $e$, $R=e\Gamma e$ is a Gorenstein toric singularity in dimension 3
(see \cite{Br,Bo2} and references therein).
Using a perfect matching $C$ on $G$, define a $\Z$-grading on $\Gamma$
by $\deg a=1$ for all $a\in C$ and $\deg a=0$ otherwise.
If both $\Gamma/(e)$ and $\Lambda=\Gamma_0$ are finite dimensional and $e\Lambda(1-e)=0$ holds,
then Theorem \ref{air theorem} gives triangle equivalences
\[\DDD^{\bo}(\mod\Lambda/(e))\simeq\underline{\CM}^{\Z}R\ \mbox{ and }\ \CCC_2(\Lambda/(e))\simeq\underline{\CM}R.\]
\end{example}

\subsection{$d$-representation-finite algebras}\label{section: d-representation finite case}

In this subsection, we study the $(d+1)$-preprojective algebras of $d$-representation-finite algebras.
We start with the following basic properties.

\begin{proposition}\cite{GLS0,I4,IO}\label{i3io}
Let $\Lambda$ be a $d$-representation-finite $k$-algebra and $\Pi=\Pi_{d+1}(\Lambda)$.
\begin{enumerate}[\rm(a)]
\item $\Pi$ is a $\Z$-graded finite dimensional self-injective $k$-algebra.
\item $\underline{\mod}^{\Z}\Pi$ has a Serre functor $(-1)[d+1]$, and $\underline{\mod}\Pi$
is $(d+1)$-Calabi-Yau.
\item $\Pi$ is a (unique) $d$-cluster tilting object in $\mod\Lambda$.
\end{enumerate}
\end{proposition}

Now we give an explicit characterization of such $\Pi$.

\begin{definition}
Let $\Gamma=\bigoplus_{i\ge0}\Gamma_i$ be a $\Z$-graded finite dimensional self-injective $k$-algebra.
We denote by $\Gamma^{\rm e}=\Gamma^{\op}\otimes_k\Gamma$ the enveloping algebra of $\Gamma$.
We say that $\Gamma$ is a \emph{stably $d$-Calabi-Yau algebra of $a$-invariant $a$}
(or \emph{Gorenstein parameter $-a$}) if
$\RHom_{\Gamma^{\rm e}}(\Gamma,\Gamma^{\rm e})(a)[d]\simeq\Gamma$ in
$\DDD_{\sg}^{\Z}(\Gamma^{\rm e})$.
\end{definition}

Now we give a homological characterization of the $(d+1)$-preprojective algebras of
$d$-representation-finite algebras as the explicit correspondence.

\begin{theorem}\cite{AO}
There exists a bijection between the set of isomorphism classes of $d$-representation-finite algebras $\Lambda$ and the set of isomorphism classes of stably $(d+1)$-Calabi-Yau self-injective algebras $\Gamma$ of $a$-invariant $-1$.
It is given by $\Lambda\mapsto\Pi_{d+1}(\Lambda)$ and $\Gamma\mapsto\Gamma_0$.
\end{theorem}

Now let $\Lambda$ be a $d$-representation-finite $k$-algebra, and $\Pi=\Pi_{d+1}(\Lambda)$.
Let $\Gamma=\underline{\End}_\Lambda(\Pi)$ be the stable $d$-Auslander algebra of $\Lambda$.
Then we have an equivalence
\begin{equation}\label{U is repetitive}
\UUU_d(\Lambda)\simeq\proj^{\Z}T(\Gamma)
\end{equation}
of additive categories. Thus we have triangle equivalences
\[\underline{\mod}^{\Z}\Pi\stackrel{\eqref{U is proj Pi}}{\simeq}\underline{\mod}
\UUU_d(\Lambda)\stackrel{\eqref{U is repetitive}}{\simeq}\underline{\mod}^{\Z}T(\Gamma)
\stackrel{\eqref{happel repetitive2}}{\simeq}\DDD^{\bo}(\mod\Gamma).\]
By Proposition \ref{i3io}(b), the automorphism $(-1)$ on $\underline{\mod}^{\Z}\Pi$ corresponds to
$\nu_{d+1}$ on $\DDD^{\bo}(\mod\Gamma)$.
Using universality of $(d+1)$-cluster categories \cite{Ke2},
we obtain a triangle equivalence $\underline{\mod}\Pi\simeq\CCC_{d+1}(\Gamma)$.
As a summary, we obtain the following.

\begin{theorem}\cite{IO}\label{io}
Let $\Lambda$ be a $d$-representation-finite $k$-algebra, $\Pi=\Pi_{d+1}(\Lambda)$,
and $\Gamma=\underline{\End}_\Lambda(\Pi)$ the stable $d$-Auslander algebra of $\Lambda$.
Then there exist triangle equivalences
\begin{eqnarray*}
\underline{\mod}^{\Z}\Pi\simeq\DDD^{\bo}(\mod\Gamma)\ \mbox{ and }\ \underline{\mod}\Pi\simeq\CCC_{d+1}(\Gamma).
\end{eqnarray*}
\end{theorem}

Applying Theorem \ref{io} for $d=1$, we obtain the following observations.

\begin{example}\cite{A,IO}
Let $\Pi$ be the preprojective algebra of a Dynkin quiver $Q$, and $\Gamma$ the stable 
Auslander algebra of $kQ$. Then there exist triangle equivalences
\begin{eqnarray*}
\underline{\mod}^{\Z}\Pi\simeq\DDD^{\bo}(\mod\Gamma)\ \mbox{ and }\ \underline{\mod}\Pi\simeq\CCC_2(\Gamma).
\end{eqnarray*}
As an application, if a quiver $Q'$ has the same underlying graph with $Q$, then the stable Auslander algebra $\Gamma'$ of $kQ'$ is derived equivalent to $\Gamma$ since $\Pi$ is common.
\end{example}

In the rest of this subsection, we discuss properties of $\Pi_{d+1}(\Lambda)$ for
a more general class of $\Lambda$.
We say that a finite dimensional $k$-algebra $\Lambda$ with $\gl\Lambda\le d$ satisfies
the \emph{vosnex property} if $\Lambda$ is $\nu_d$-finite and satisfies
$\Hom_{\DDD^{\bo}(\mod\Lambda)}(\UUU_d(\Lambda)[i],\UUU_d(\Lambda))=0$ for all $1\le i\le d-2$.
This is automatic if $d=1,2$ or $\Lambda$ is $d$-representation-finite.
In this case, the following generalization of Theorem \ref{io} holds.

\begin{theorem}\cite{IO}\label{io2}
Let $\Lambda$ be a finite dimensional $k$-algebra with $\gl\Lambda\le d$ satisfying
the vosnex property.
Then $\Pi=\Pi_{d+1}(\Lambda)$ is 1-Iwanaga-Gorenstein, $\Gamma=\underline{\End}_{\Lambda}(\Pi)$
satisfies $\gl\Gamma\le d+1$, and there exist triangle equivalences
\begin{eqnarray*}
\underline{\CM}^{\Z}\Pi\simeq\DDD^{\bo}(\mod\Gamma)\ \mbox{ and }\ \underline{\CM}\Pi\simeq\CCC_{d+1}(\Gamma).
\end{eqnarray*}
\end{theorem}

For more general $\Lambda$, we refer to \cite{Be} for some properties of $\Pi_{d+1}(\Lambda)$.

\subsection{Preprojective algebras and Coxeter groups}\label{section: Truncated preprojective algebras}
We discuss a family of finite dimensional $k$-algebras constructed from preprojective algebras and
Coxeter groups.

Let $Q$ be an acyclic quiver and $\Pi$ the preprojective algebra of $kQ$.
The \emph{Coxeter group} of $Q$ is generated by $s_i$ with $i\in Q_0$, and the relations are the following.
\begin{itemize}
\item $s_i^2=1$ for all $i\in Q_0$.
\item $s_is_j=s_js_i$ if there is no arrow between $i$ and $j$ in $Q$.
\item $s_is_js_i=s_js_is_j$ if there is precisely one arrow between $i$ and $j$ in $Q$.
\end{itemize}
Let $w\in W$. An expression $w=s_{i_1}s_{i_2}\cdots s_{i_\ell}$ of $w$ is called \emph{reduced} if
$\ell$ is minimal among all expressions of $w$.
For $i\in Q_0$, let $I_i$ be the two-sided ideal of $\Pi$ generated by the idempotent $1-e_i$.
For a reduced expression $w=s_{i_1}\cdots s_{i_\ell}$, we define a two-sided ideal of $\Pi$ by
\[I_w:=I_{i_1}I_{i_2}\cdots I_{i_\ell}.\]
This is independent of the choice of the reduced expression of $w$.
The corresponding factor algebra $\Pi_w:=\Pi/I_w$ is a finite dimensional $k$-algebra.
It enjoys the following remarkable properties.

\begin{theorem}\cite{BIRS,GLS1,ART}\label{air}
\begin{enumerate}[\rm(a)]
\item $\Pi_w$ is a $1$-Iwanaga-Gorenstein algebra.
\item $\underline{\CM}\Pi_w$ is a 2-Calabi-Yau triangulated category.
\item There exists a 2-cluster tilting object $\bigoplus_{j=1}^\ell e_{i_j}\Pi_{s_{i_j}\cdots s_{i_\ell}}$
in $\CM\Pi_w$.
\item There exists a triangle equivalence $\underline{\CM}\Pi_w\simeq\CCC_2(\Lambda)$ for
some algebra $\Lambda$.
\end{enumerate}
\end{theorem}

Therefore it is natural to expect that there exists a triangle equivalence
$\underline{\CM}^{\Z}\Pi_w\simeq\DDD^{\bo}(\mod\Lambda')$ for some algebra $\Lambda'$.
In fact, the following results are known, where we refer to \cite{Ki1,Ki2} for the definitions of
\emph{$c$-sortable}, \emph{$c$-starting} and \emph{$c$-ending}.

\begin{theorem}\label{ki}
Let $w=s_{i_1}\cdots s_{i_\ell}$ be a reduced expression of $w\in W$.
\begin{enumerate}[\rm(a)]
\item \cite{Ki1} If $w$ is $c$-sortable, then $\underline{\CM}^{\Z}\Pi_w$ has a tilting object
$\bigoplus_{i>0}\Pi_w(i)_{\ge0}$.
\item \cite{Ki2} $\underline{\CM}^{\Z}\Pi_w$ has a silting object $\bigoplus_{j=1}^\ell e_{i_j}\Pi_{s_{i_j}\cdots s_{i_\ell}}$.
This is a tilting object if the reduced expression is $c$-starting or $c$-ending.
\end{enumerate}
\end{theorem}

We end this section with posing the following natural question on `higher cluster combinatorics' (e.g.\ \cite{OT}), which will be related to derived equivalences of Calabi-Yau algebras since our $I_w$ is a tilting object in $\KKK^{\bo}(\proj\Pi)$ if $Q$ is non-Dynkin.

\begin{problem}
Are there similar results to Theorems \ref{air} and \ref{ki} for higher preprojective algebras?
What kind of combinatorial structure will appear instead of the Coxeter groups?
\end{problem}

\section{Geigle-Lenzing complete intersections}

Weighted projective lines of Geigle-Lenzing \cite{GL1} are one of the basic objects in representation theory.
For example, the simplest class of weighted projective lines gives us simple singularities
in dimension 2 as certain Veronese subrings.
We introduce a higher dimensional generalization of weighted projective lines following \cite{HIMO}.

\subsection{Basic properties}
For a field $k$ and an integer $d\ge1$, we consider a polynomial algebra $C=k[T_0,\ldots,T_d]$.
For $n\ge0$, let $\ell_1,\ldots,\ell_n$ be linear forms in $C$ and $p_1,\ldots,p_n$ 
positive integers. For simplicity, we assume $p_i\ge2$ for all $i$. Let
\begin{equation*}\label{show R}
R=C[X_1, \ldots, X_n]/(X_i^{p_i}-\ell_i\mid 1\le i\le n)
\end{equation*}
be the factor algebra of the polynomial algebra $C[X_1,\ldots,X_n]$, and
\begin{equation*}\label{show L}
\L = \langle \x_1,\ldots,\x_n, \c \rangle / \langle p_i \x_i - \c \mid 1\le i \le n\rangle.
\end{equation*}
the factor group of the free abelian group $\langle \x_1,\ldots,\x_n, \c \rangle$.
Then $\L$ is an abelian group of rank 1 with torsion elements in general, and
$R$ is $\L$-graded by $\deg T_i=\c$ for all $0\le i\le d$ and $\deg X_i=\x_i$ for all $1\le i\le n$.

We call the pair $(R,\L)$ a \emph{Geigle-Lenzing} (\emph{GL}) \emph{complete intersection}
if $\ell_1,\ldots,\ell_n$ are in general position in the sense that each set of at most $d+1$ elements
from $\ell_1,\ldots,\ell_n$ is linearly independent.
We give some basic properties.

\begin{proposition}
Let $(R,\L)$ be a GL complete intersection.
\begin{enumerate}[\rm(a)]
\item $X_1^{p_1}-\ell_1,\ldots,X_n^{p_n}-\ell_n$ is a $C[X_1,\ldots,X_n]$-regular sequence.
\item $R$ is a complete intersection ring with $\dim R=d+1$ and has an $a$-invariant
\[\w=(n-d-1)\c-\sum_{i=1}^n\x_i.\]
\item After a suitable linear transformation of variables $T_0,\ldots,T_d$, we have
\begin{equation*}\label{normalized form}
R = \begin{cases} k[X_1, \ldots, X_n, T_n, \ldots, T_d] & \text{ if } n \leq d+1, \\
k[X_1, \ldots, X_n] / (X_i^{p_i} - \sum_{j=1}^{d+1}\lambda_{i,j-1}X_j^{p_j}) \mid d+2 \leq i \leq n) & \text{ if } n \geq d+2. \end{cases}
\end{equation*}
\item $R$ is regular if and only if $R$ is a polynomial algebra if and only if $n\le d+1$.
\item $\CM^{\L}R=\CM^{\L}_0R$ holds, and $\underline{\CM}^{\L}R$ has a Serre functor $(\w)[d]$ (Theorem \ref{AR duality for CM}).
\end{enumerate}
\end{proposition}

Let $\delta\colon \L\to\Q$ be a group homomorphism given by $\delta(\x_i)=\frac{1}{p_i}$ and $\delta(\c)=1$.
We consider the following trichotomy given by the sign of
$\delta(\w)=n-d-1-\sum_{i=1}^n\frac{1}{p_i}$. For example, $(R,\L)$ is Fano if $n\le d+1$.
\[\begin{array}{|c||c|c|c|}\hline
\delta(\w)&<0&=0&>0\\ \hline\hline
(R,\L)&\mbox{Fano}&\mbox{Calabi-Yau}&\mbox{anti-Fano}\\ \hline
d=1&\mbox{domestic}&\mbox{tubular}&\mbox{wild}\\ \hline
\end{array}\]
In the classical case $d=1$, the ring $R$ has been studied in the context of
weighted projective lines. The above trichotomy is given explicitly as follows.
\begin{enumerate}[$\bullet$]
\item 5 types for domestic: $n\le 2$, $(2,2,p)$, $(2,3,3)$, $(2,3,4)$ and $(2,3,5)$.
\item 4 types for tubular: $(3,3,3)$, $(2,4,4)$, $(2,3,6)$ and $(2,2,2,2)$.
\item All other types are wild.
\end{enumerate}
There is a close connection between domestic type and simple singularities.
The following explains Corollary \ref{simple surface singularity}, where 
$R^{(\w)}=\bigoplus_{i\in\Z}R_{i\w}$ is the Veronese subring.

\begin{theorem}\cite{GL2}\label{R^w is simple}
If $(R,\L)$ is domestic, then $R^{(\w)}$ is a simple singularity $k[x,y,z]/(f)$ in dimension 2,
and we have an equivalence $\CM^{\L}R\simeq\CM^{\Z}R^{(\w)}$. The AR quiver is
$\Z Q$, where $Q$ is given by the following table.
\[\begin{array}{|c||c|c|c|c|c|}\hline 
(p_1,\ldots,p_n)&x&y&z&f&Q\\ \hline\hline
(p,q)&X_1X_2&X_2^{p+q}&X_1^{p+q}&x^{p+q}-yz&\widetilde{\A}_{p,q}\\ \hline
(2,2,2p)&X_3^2&X_1^2&X_1X_2X_3&x(y^2+x^py)+z^2&\widetilde{\D}_{2p+2}\\ \hline
(2,2,2p+1)&X_3^2&X_1X_2&X_1^2X_3&x(y^2+x^pz)+z^2&\widetilde{\D}_{2p+3}\\ \hline
(2,3,3)&X_1&X_2X_3&X_2^3&x^2z+y^3+z^2&\widetilde{\E}_6\\ \hline
(2,3,4)&X_2&X_3^2&X_1X_3&x^3y+y^3+z^2&\widetilde{\E}_7\\ \hline
(2,3,5)&X_3&X_2&X_1&x^5+y^3+z^2&\widetilde{\E}_8\\ \hline
\end{array}\]
\end{theorem}

\subsection{Cohen-Macaulay representations}

To study the category $\CM^{\L}R$, certain finite dimensional algebras
play an important role. For a finite subset $I$ of $\L$, let
\[A^I=\bigoplus_{\x,\y\in I}R_{\x-\y}.\]
We define the multiplication in $A^I$ by
$(r_{\x,\y})_{\x,\y\in I}\cdot(r'_{\x,\y})_{\x,\y\in I}=(\sum_{\z\in I}r_{\x,\z}r'_{\z,\y})_{\x,\y\in I}$.
Then $A^I$ forms a finite dimensional $k$-algebra called the \emph{$I$-canonical algebra}.

We define a partial order $\le$ on $\L$ by writing $\x\le\y$ if $\y-\x$ belongs to $\L_+$,
where $\L_+$ is the submonoid of $\L$ generated by $\c$ and $\x_i$ for all $i$. 
For $\x\in\L$, let $[0,\x]$ be the interval in $\L$, and $A^{[0,\x]}$ the $[0,\x]$-canonical algebra.
We call
\[A^{\rm CM}=A^{[0,d\c+2\w]}\]
the \emph{CM-canonical algebra}.

\begin{example}\label{example of canonical}
The equality $d\c+2\w=(n-d-2)\c+\sum_{i=1}^n(p_i-2)\x_i$ holds.
\begin{enumerate}[\rm(a)]
\item If $n\le d+1$, then $A^{\rm CM}=0$. If $n=d+2$, then $A^{\rm CM}=\bigotimes_{i=1}^{n+2}k\A_{p_i-1}$.
\item If $n=d+3$ and $p_i=2$ for all $i$, then $A^{\rm CM}$ has the left quiver below.
\item If $d=1$, $n=4$ and $(p_i)_{i=1}^4=(2,2,2,3)$, then $A^{\rm CM}$ has the right quiver below.
{\small
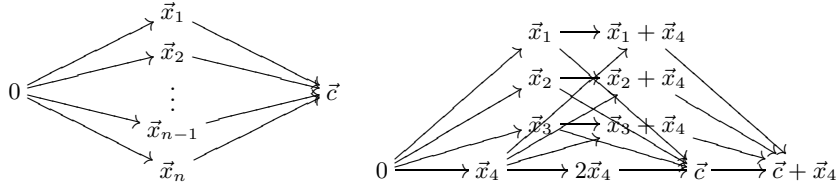
\[\begin{xy} 0;<4pt,0pt>:<0pt,3pt>:: 
(0,10) *+{0} ="0",
(15,0) *+{\x_n} ="n",
(15,5) *+{\x_{n-1}} ="n-1",
(15,10) *+{\vdots},
(15,15) *+{\x_2} ="2",
(15,20) *+{\x_1}="1",
(30,10) *+{\c}="c",
"0", {\ar"1"},
"0", {\ar"2"},
"0", {\ar"n-1"},
"0", {\ar"n"},
"1", {\ar"c"},
"2", {\ar"c"},
"n-1", {\ar"c"},
"n", {\ar"c"},
\end{xy}\ \ \ 
\begin{xy} 0;<4pt,0pt>:<0pt,3.5pt>:: 
(-20,0) *+{0} ="0",
(-10,0) *+{\x_4} ="4",
(-5,15) *+{\x_1} ="1",
(-5,10) *+{\x_2} ="2",
(-5,5) *+{\x_3} ="3",
(0,0) *+{2\x_4} ="44",
(5,15) *+{\x_1+\x_4} ="14",
(5,10) *+{\x_2+\x_4} ="24",
(5,5) *+{\x_3+\x_4} ="34",
(10,0) *+{\c} ="c",
(20,0) *+{\c+\x_4} ="c4",
"0", {\ar"1"},
"0", {\ar"2"},
"0", {\ar"3"},
"0", {\ar"4"},
"1", {\ar"14"},
"1", {\ar"c"},
"2", {\ar"24"},
"2", {\ar"c"},
"3", {\ar"34"},
"3", {\ar"c"},
"4", {\ar"14"},
"4", {\ar"24"},
"4", {\ar"34"},
"4", {\ar"44"},
"44", {\ar"c"},
"c", {\ar"c4"},
"14", {\ar"c4"},
"24", {\ar"c4"},
"34", {\ar"c4"},
\end{xy}\]}
\end{enumerate}
\end{example}

The following is a main result in this section.
\begin{theorem}
Let $(R,\L)$ be a GL complete intersection. Then there is a triangle equivalence
\[\underline{\CM}^{\L}R\simeq\DDD^{\bo}(\mod A^{\rm CM}).\]
In particular, $\underline{\CM}^{\L}R$ has a tilting object.
\end{theorem}

The case $n=d+2$ was shown in \cite{KLM} ($d=1$) and \cite{FU}.
An important tool in the proof is an $\L$-analogue of Theorem \ref{orlov embedding}.

As an application, one can immediately obtain the following analogue of
Theorem \ref{bgs} by using the knowledge on $A^{\rm CM}$ in representation theory,
where we call $(R,\L)$ \emph{CM-finite} if there are only finitely many isomorphism classes
of indecomposable objects in $\CM^{\L}R$ up to degree shift (cf.\ Definition \ref{define CM-finite}).

\begin{corollary}
Let $(R,\L)$ be a GL complete intersection.
Then $(R,\L)$ is CM-finite if and only if one of the following conditions hold.  
\begin{enumerate}[\rm(i)]
\item $n\le d+1$.
\item $n=d+2$, and $(p_1,\ldots,p_n)=(2,\ldots,2,p_n)$, $(2,\ldots,2,3,3)$, $(2,\ldots,2,3,4)$
or $(2,\ldots,2,3,5)$ up to permutation.
\end{enumerate}
\end{corollary}

We call a GL complete intersection $(R,\L)$
\emph{$d$-CM-finite} if there exists a $d$-cluster tilting subcategory $\CC$ of $\CM^{\L}R$ such that there
are only finitely many isomorphism classes of indecomposable objects in $\CC$ up to degree shift
(cf.\ Definition \ref{define d-CM-finite}).
Now we discuss which GL complete intersections are $d$-CM-finite.
Our Theorem \ref{C_d has d-CT} gives the following sufficient condition, where a tilting object is called \emph{$d$-tilting}
if the endomorphism algebra has global dimension at most $d$.

\begin{proposition}\label{d-tilting object}
If $\underline{\CM}^{\L}R$ has a $d$-tilting object $U$, then $(R,\L)$ is $d$-CM-finite and $\CM^{\L}R$
has the $d$-cluster tilting subcategory $\add\{U(\ell\w), R(\x)\mid\ell\in\Z, \x\in\L\}$.
\end{proposition}

Therefore the following problem is of our interest.

\begin{problem}
When does $\underline{\CM}^{\L}R$ have a $d$-tilting object?
Equivalently, when is $A^{\rm CM}$ derived equivalent to an algebra $\Lambda$ with $\gl\Lambda\le d$?
\end{problem}

Applying Tate's DG algebra resolutions \cite{T}, we can calculate $\gl A^{\rm CM}$.
Note that any element $\x\in\L$ can be written uniquely as $\x=a\c+\sum_{i=1}^na_i\x_i$
for $a\in\Z$ and $0\le a_i\le p_i-1$, which is called the \emph{normal form} of $\x$.

\begin{theorem}\label{calculate gl.dim}
\begin{enumerate}[\rm(a)]
\item Write $\x\in\L_+$ in normal form $\x=a\c+\sum_{i=1}^na_i\x_i$. Then
\[\gl A^{[0,\x]}=\begin{cases}
\min\{d+1,a+\#\{i\mid a_i\neq0\}\}&\mbox{ if }\ n\le d+1,\\
2a+\#\{i\mid a_i\neq0\}&\mbox{ if }\ n\ge d+2.\end{cases}\]
\item If $n\ge d+2$, then $A^{\rm CM}$ has global dimension $2(n-d-2)+\#\{i\mid p_i\ge 3\}$.
\end{enumerate}
\end{theorem}

We obtain the following examples from Theorem \ref{calculate gl.dim} and the fact that
$k\A_2\otimes_kk\A_m$ is derived equivalent to $k\D_4$ if $m=2$, $k\E_6$ if $m=3$, and $k\E_8$ if $m=4$.

\begin{example}\label{d-tilting exist}
In the following cases, $\underline{\CM}^{\L}R$ has a $d$-tilting object.
\begin{enumerate}[\rm(i)]
\item $n\le d+1$.
\item $n=d+2\ge3$ and $(p_1,p_2,p_3)=(2,2,p_3)$, $(2,3,3)$, $(2,3,4)$ or $(2,3,5)$.
\item $n=d+2\ge4$ and $(p_1,p_2,p_3,p_4)=(3,3,p_3,p_4)$ with $p_3,p_4\in\{3,4,5\}$.
\item $\#\{i\mid p_i=2\}\ge3(n-d)-4$.
\end{enumerate}
\end{example}

The following gives a necessary condition for the existence of $d$-tilting object.

\begin{proposition}\label{d-tilting implies fano}
If $\underline{\CM}^{\L}R$ has a $d$-tilting object, then $(R,\L)$ is Fano.
\end{proposition}

Note that the converse is not true. For example, let $d=2$ and $(2,5,5,5)$.
Then $(R,\L)$ is Fano since $\delta(\w)=-\frac{1}{10}$. On the other hand,
$A^{\rm CM}=\bigotimes_{i=1}^3k\A_4$
satisfies $\nu^5=[9]$. One can show that $A^{\rm CM}$ is not derived equivalent
to an algebra $\Lambda$ with $\gl\Lambda\le 2$ by using the inequality $2(5-1)<9$.

\subsection{Geigle-Lenzing projective spaces}

Let $(R,\L)$ be a GL complete intersection. Recall that $\mod^{\L}_0R$ is the Serre subcategory
of $\mod^{\L}R$ consisting of finite dimensional modules. We consider the quotient category
\[\coh\X=\qgr R=\mod^{\L}R/\mod^{\L}_0R.\]
We call objects in $\coh\X$ \emph{coherent sheaves} on the \emph{GL projective space} $\X$.
We can regard $\X$ as the quotient stack $[(\Spec R\setminus\{R_+\})/\Spec k[\L]]$
for $R_+=\bigoplus_{\x>0}R_{\x}$. For example, if $n=0$, then $\X$ is the projective space $\P^d$.

We study the bounded derived category $\DDD^{\bo}(\coh\X)$, which is canonically triangle equivalent to
the Verdier quotient $\DDD^{\bo}(\mod^{\L}R)/\DDD^{\bo}(\mod_0^{\L}R)$.
The duality $(-)^*=\RHom_R(-,R)\colon \DDD^{\bo}(\mod^{\L}R)\to\DDD^{\bo}(\mod^{\L}R)$
induces a duality $(-)^*\colon \DDD^{\bo}(\coh\X)\to\DDD^{\bo}(\coh\X)$.
We define the category of \emph{vector bundles} on $\X$ as
\[\vect\X=\coh\X\cap(\coh\X)^*.\]
The composition $\CM^{\L}R\subset\mod^{\L}R\to\coh\X$ is fully faithful, and we can regard
$\CM^{\L}R$ as a full subcategory of $\vect\X$.
We have $\CM^{\L}R=\vect\X$ if $d=1$, but this is not the case if $d\ge 2$. In fact, we have equalities
\begin{eqnarray}\label{CM and vect}
\CM^{\L}R&=&\{X\in\vect\X\mid\forall\x\in\L,\ 1\le i\le d-1,\ \Ext^i_{\X}(\OO(\x),X)=0\}\\ \notag
&=&\{X\in\vect\X\mid\forall\x\in\L,\ 1\le i\le d-1,\ \Ext^i_{\X}(X,\OO(\x))=0\},
\end{eqnarray}
where $\OO(\x)=R(\x)$. Now we define the \emph{$d$-canonical algebra} by
\[A^{\rm ca}=A^{[0,d\c]}.\]

\begin{example}
\begin{enumerate}[\rm(a)]
\item If $d=1$, then $A^{\rm ca}$ is precisely the canonical algebra of Ringel \cite{R}.
It is given by the following quiver with relations $x_i^{p_i}=\lambda_{i0}x_1^{p_1}+\lambda_{i1}x_2^{p_2}$ for any $i$ with $3\le i\le n$.
{\small
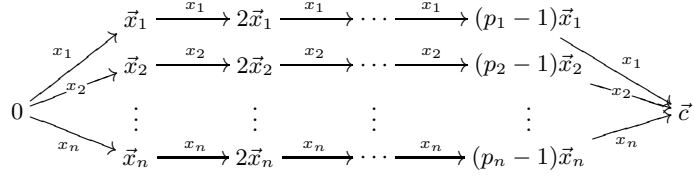
\[
\begin{xy} 0;<4.5pt,0pt>:<0pt,3.5pt>:: 
(0,0) *+{0} ="0",
(10,10) *+{\x_1} ="1",
(10,5) *+{\x_2} ="2",
(10,0) *+{\vdots} ="e",
(10,-5) *+{\x_n} ="n",
(20,10) *+{2\x_1} ="11",
(20,5) *+{2\x_2} ="22",
(20,0) *+{\vdots} ="ee",
(20,-5) *+{2\x_n} ="nn",
(30,10) *+{\cdots} ="1e1",
(30,5) *+{\cdots} ="2e2",
(30,0) *+{\vdots} ="eee",
(30,-5) *+{\cdots} ="nen",
(43,10) *+{(p_1-1)\x_1} ="p1",
(43,5) *+{(p_2-1)\x_2} ="p2",
(43,0) *+{\vdots} ="pe",
(43,-5) *+{(p_n-1)\x_n} ="pn",
(56,0) *+{\c} ="c",
"0", {\ar^{x_1}"1"},
"0", {\ar|{x_2}"2"},
"0", {\ar_{x_n}"n"},
"1", {\ar^{x_1}"11"},
"2", {\ar^{x_2}"22"},
"n", {\ar^{x_n}"nn"},
"11", {\ar^{x_1}"1e1"},
"22", {\ar^{x_2}"2e2"},
"nn", {\ar^{x_n}"nen"},
"1e1", {\ar^{x_1\rule{15pt}{0pt}}"p1"},
"2e2", {\ar^{x_2\rule{15pt}{0pt}}"p2"},
"nen", {\ar^{x_n\rule{15pt}{0pt}}"pn"},
"p1", {\ar^(.6){x_1}"c"},
"p2", {\ar|(.6){x_2}"c"},
"pn", {\ar_(.6){x_n}"c"},
\end{xy}
\]}
\item If $n=0$, then $A^{\rm ca}$ is the Beilinson algebra.
\item If $d=2$, $n=3$ and $(p_i)_{i=1}^3=(2,2,2)$, then $A^{\rm ca}$ has the left quiver below.
\item If $d=2$, $n=4$ and $(p_i)_{i=1}^4=(2,2,2,2)$, then $A^{\rm ca}$ has the right quiver below.
{\tiny
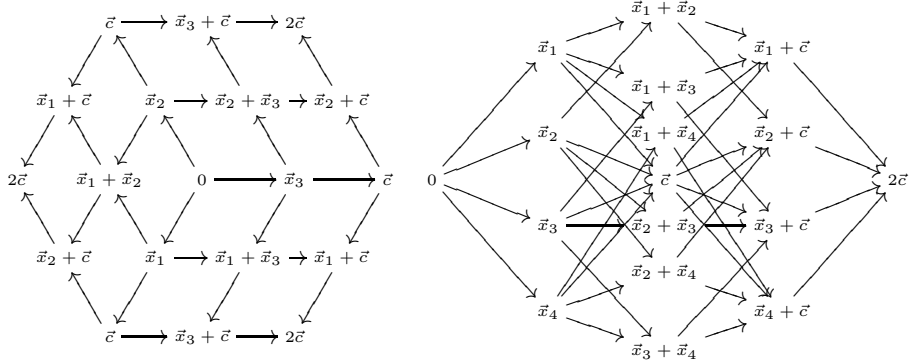
\[
\begin{xy}
0;<35pt,0cm>:<-17.5pt,29.75pt>:: 
(0,0) *+{0} ="000",
(1,0) *+{\x_3} ="100",
(2,0) *+{\c} ="200",
(0,1) *+{\x_2} ="010",
(0,2) *+{\c} ="020",
(-1,-1) *+{\x_1} ="001",
(-2,-2) *+{\c} ="002",
(1,1) *+{\x_2+\x_3} ="110",
(2,1) *+{\x_2+\c} ="210",
(1,2) *+{\x_3+\c} ="120",
(2,2) *+{2\c} ="220",
(0,-1) *+{\x_1+\x_3} ="101",
(1,-1) *+{\x_1+\c} ="201",
(-1,-2) *+{\x_3+\c} ="102",
(0,-2) *+{2\c} ="202",
(-1,0) *+{\x_1+\x_2} ="011",
(-1,1) *+{\x_1+\c} ="021",
(-2,-1) *+{\x_2+\c} ="012",
(-2,0) *+{2\c} ="022",
"000", {\ar"100"},
"100", {\ar"200"},
"010", {\ar"110"},
"110", {\ar"210"},
"020", {\ar"120"},
"120", {\ar"220"},
"000", {\ar"010"},
"010", {\ar"020"},
"100", {\ar"110"},
"110", {\ar"120"},
"200", {\ar"210"},
"210", {\ar"220"},
"001", {\ar"101"},
"101", {\ar"201"},
"002", {\ar"102"},
"102", {\ar"202"},
"000", {\ar"001"},
"001", {\ar"002"},
"100", {\ar"101"},
"101", {\ar"102"},
"200", {\ar"201"},
"201", {\ar"202"},
"001", {\ar"011"},
"011", {\ar"021"},
"002", {\ar"012"},
"012", {\ar"022"},
"010", {\ar"011"},
"011", {\ar"012"},
"020", {\ar"021"},
"021", {\ar"022"},
\end{xy}\ \ \ 
\begin{xy} 0;<2.2pt,0pt>:<0pt,2.5pt>:: 
(-40,0) *+{0} ="0",
(-20,20) *+{\x_1} ="1",
(-20,7) *+{\x_2} ="2",
(-20,-7) *+{\x_3} ="3",
(-20,-20) *+{\x_4} ="4",
(0,26) *+{\x_1+\x_2} ="12",
(0,14) *+{\x_1+\x_3} ="13",
(0,7) *+{\x_1+\x_4} ="14",
(0,-7) *+{\x_2+\x_3} ="23",
(0,-14) *+{\x_2+\x_4} ="24",
(0,-26) *+{\x_3+\x_4} ="34",
(0,0) *+{\c} ="c",
(20,20) *+{\x_1+\c} ="c1",
(20,7) *+{\x_2+\c} ="c2",
(20,-7) *+{\x_3+\c} ="c3",
(20,-20) *+{\x_4+\c} ="c4",
(40,0) *+{2\c} ="2c",
"0", {\ar"1"},
"0", {\ar"2"},
"0", {\ar"3"},
"0", {\ar"4"},
"1", {\ar"c"},
"1", {\ar"12"},
"1", {\ar"13"},
"1", {\ar"14"},
"2", {\ar"12"},
"2", {\ar"c"},
"2", {\ar"23"},
"2", {\ar"24"},
"3", {\ar"13"},
"3", {\ar"23"},
"3", {\ar"c"},
"3", {\ar"34"},
"4", {\ar"14"},
"4", {\ar"24"},
"4", {\ar"34"},
"4", {\ar"c"},
"12", {\ar"c2"},
"12", {\ar"c1"},
"13", {\ar"c3"},
"13", {\ar"c1"},
"14", {\ar"c4"},
"14", {\ar"c1"},
"23", {\ar"c3"},
"23", {\ar"c2"},
"24", {\ar"c4"},
"24", {\ar"c2"},
"34", {\ar"c4"},
"34", {\ar"c3"},
"c", {\ar"c1"},
"c", {\ar"c2"},
"c", {\ar"c3"},
"c", {\ar"c4"},
"c1", {\ar"2c"},
"c2", {\ar"2c"},
"c3", {\ar"2c"},
"c4", {\ar"2c"},
\end{xy}
\]}
\end{enumerate}
\end{example}

As in the case of $\underline{\CM}^{\L}R$ and $A^{\rm CM}$, we obtain the following results.

\begin{theorem}
Let $\X$ be a GL projective space.
Then there is a triangle equivalence
\[\DDD^{\bo}(\coh\X)\simeq\DDD^{\bo}(\mod A^{\rm ca}).\]
Moreover $\DDD^{\bo}(\coh\X)$ has a tilting bundle $\bigoplus_{\x\in[0,d\c]}\OO(\x)$.
\end{theorem}

Some cases were known before ($n=0$ \cite{Bei}, $d=1$ \cite{GL1}, $n\le d+1$ \cite{Ba},
$n=d+2$ \cite{IU}).
An important tool in the proof is again an $\L$-analogue of Theorem \ref{orlov embedding}.

We call $\X$ \emph{vector bundle finite}  (\emph{VB-finite}) if there are only finitely many isomorphism
classes of indecomposable objects in $\vect\X$ up to degree shift.
There is a complete classification: $\X$ is VB-finite if and only if $d=1$ and $\X$ is domestic.

We call $\X$ \emph{$d$-VB-finite} if there exists a $d$-cluster tilting subcategory $\CC$ of $\vect\X$
such that there are only finitely many isomorphism classes of indecomposable objects in $\CC$
up to degree shift. In the rest, we discuss which GL projective spaces are $d$-VB-finite.
We start with the following relationship between $d$-cluster tilting subcategories of
$\CM^{\L}R$ and $\vect\X$, which follows from \eqref{CM and vect}.

\begin{proposition}\label{d-CM-finite implies d-VB-finite}
The $d$-cluster tilting subcategories of $\CM^{\L}R$ are precisely the $d$-cluster tilting subcategories of $\vect\X$ containing $\OO(\x)$ for all $\x\in\L$.
Therefore, if $(R,\L)$ is $d$-CM-finite, then $\X$ is $d$-VB-finite.
\end{proposition}

For example, if $n\le d+1$, then $\CM^{\L}R=\proj^{\L}R$ is a $d$-cluster tilting subcategory of itself, and hence $\vect\X$ has a $d$-cluster tilting subcategory $\add\{\OO(\x)\mid\x\in\L\}$.
This implies Horrocks' splitting criterion for $\vect\P^d$ \cite{OSS}.

We give another sufficient condition for $d$-VB-finiteness. Recall that we call a tilting object $V$
in $\DDD^{\bo}(\coh\X)$ \emph{$d$-tilting} if $\gl\End_{\DDD^{\bo}(\coh\X)}(V)\le d$.

\begin{proposition}\label{d-tilting sheaf}
Let $\X$ be a GL projective space, and $V$ a $d$-tilting object in $\DDD^{\bo}(\coh\X)$.
\begin{enumerate}[\rm(a)]
\item (cf.\ Example \ref{example of d-hereditary}(b)) $\gl\End_{\DDD^{\bo}(\coh\X)}(V)=d$ holds.
If $V\in\coh\X$, then $\End_{\X}(T)$ is a $d$-representation-infinite algebra.
\item If $V\in\vect\X$, then $\X$ is $d$-VB-finite and $\vect\X$
has the $d$-cluster tilting subcategory $\add\{V(\ell\w)\mid\ell\in\Z\}$.
\end{enumerate}
\end{proposition}

Therefore it is natural to ask when $\X$ has a $d$-tilting bundle, or equivalently, when $A^{\rm ca}$ is
derived equivalent to an algebra $\Lambda$ with $\gl\Lambda=d$.
It follows from Theorem \ref{calculate gl.dim}(a) that
\[\gl A^{\rm ca}=\begin{cases}
d&\mbox{ if }\ n\le d+1,\\
2d&\mbox{ if }\ n\ge d+2.\end{cases}\]
Thus, if $n\le d+1$, then $\X$ has a $d$-tilting bundle.
Using Example \ref{d-tilting exist} and some general results on matrix factorizations, we have more examples.

\begin{theorem}
In the following cases, $\X$ has a $d$-tilting bundle.
\begin{enumerate}[\rm(i)]
\item $n\le d+1$.
\item $n=d+2\ge3$ and $(p_1,p_2,p_3)=(2,2,p_3)$, $(2,3,3)$, $(2,3,4)$ or $(2,3,5)$.
\item $n=d+2\ge4$ and $(p_1,p_2,p_3,p_4)=(3,3,p_3,p_4)$ with $p_3,p_4\in\{3,4,5\}$.
\end{enumerate}
\end{theorem}

As in the previous subsection, we have the following necessary condition.

\begin{proposition}\label{d-tilting implies fano2}
If $\X$ has a $d$-tilting bundle, then $\X$ is Fano.
\end{proposition}

Some of our results in this section can be summarized as follows.
\[\xymatrix@C=4em@R=1em{
{\begin{array}{c}\mbox{$(R,\L)$ is}\\ \mbox{$d$-CM-finite}\end{array}}\ar@{=>}[d]_{\rm Prop. \ref{d-CM-finite implies d-VB-finite}}&
{\begin{array}{c}\mbox{$\underline{\CM}^{\L}R$ has a}\\ \mbox{$d$-tilting object}\end{array}}\ar@{=>}[l]_{\rm Prop. \ref{d-tilting object}}
\ar@{=>}[r]^{\rm Prop. \ref{d-tilting implies fano}}&\mbox{Fano}\\
{\begin{array}{c}\mbox{$\X$ is}\\ \mbox{$d$-VB-finite}\end{array}}&
{\begin{array}{c}\mbox{$\X$ has a}\\ \mbox{$d$-tilting bundle}\end{array}}
\ar@{=>}[l]_(.55){\rm Prop. \ref{d-tilting sheaf}}
\ar@{=>}[r]^(.4){\rm Prop. \ref{d-tilting sheaf}}
\ar@{=>}[ru]^{\rm Prop. \ref{d-tilting implies fano2}}
&{\begin{array}{c}\mbox{$\X$ is derived equivalent}\\
\mbox{to a $d$-representation}\\
\mbox{infinite algebra}
\end{array}}
}\]
It is important to understand the precise relationship between these conditions.
We refer to \cite{C,BHI} for results on existence of $d$-tilting bundles on more general varieties and stacks.

\end{document}